\documentclass[12pt]{article}

\usepackage[ansinew]{inputenc}

\usepackage{amsmath} \usepackage{amsthm}
\usepackage{amsfonts}\usepackage{amssymb}
\usepackage{graphicx}
\usepackage{color}

\newtheorem{theorem}{Theorem} \newtheorem{prop}[theorem]{Proposition}
\newtheorem{lemma}[theorem]{Lemma} 
\newtheorem{cor}[theorem]{Corollary} \newtheorem{rmk}[theorem]{Remark}

\newcommand{\eps}{\varepsilon} 

\newcommand{\indic}[1]{\mathbf{1}_{\{#1\}}}

 \DeclareMathOperator\cov{{cov}}

\def\mn{\medskip\noindent}

\def\beq{\begin{equation}}
\def\eeq{\end{equation}}
\def\beqa{\begin{eqnarray}}
\def\eeqa{\end{eqnarray}}
\def\beqax{\begin{eqnarray*}}
\def\eeqax{\end{eqnarray*}}
\def\sqz{\kern -0.2em}

\def\ep{\varepsilon}

\def\si{\sigma}

\def\om{{\omega}}
\def\Om{\Omega}
\def\Up{\Upsilon}

\def\cA{\mathcal{A}}
\def\cB{\mathcal{B}}
\def\cC{\mathcal{C}}
\def\cD{\mathcal{D}}
\def\cE{\mathcal{E}}
\def\cG{\mathcal{G}}
\def\cF{\mathcal{F}}
\def\cH{\mathcal{H}}
\def\cM{\mathcal{M}}

\def\cR{\mathcal{R}}

\def\cL{\mathcal{L}}

\def\limt{\lim_{t\to\infty}}
\def\limn{\lim_{n\to\infty}}

\def\wt{\widetilde}

\def\sobre#1#2{\lower 1ex \hbox{ $#1 \atop #2 $ } }
\def\bajo#1#2{\raise 1ex \hbox{ $#1 \atop #2 $ } }

\def\eps{\varepsilon}

\newcommand{\bm}{\mathbb{W}}
\newcommand{\be}{\mathbb{Y}}

\def\R{\mathbb{R}}

\begin{document}

\title{Random paths with bounded local time}

\author{Itai Benjamini$^1$ and Nathana\"el Berestycki$^2$}

\date{April 21, 2010}
\maketitle

\centerline{\textbf{Abstract.}}

\mn

\mn We consider one-dimensional Brownian motion conditioned (in a suitable sense) to have a local time at every point and at every moment bounded by some fixed constant. Our main result shows that a phenomenon of entropic repulsion occurs: that is, this process is ballistic and has an asymptotic velocity approximately 4.58$\ldots$ as high as required by the conditioning (the exact value of this constant involves the first zero of a Bessel function). We also study the random walk case and show that the process is asymptotically ballistic but with an unknown speed.

\vfill

\mn 1.
 \texttt{itai.benjamini@weizmann.ac.il}. Weizmann Institute of
 Science. Rehovot, Israel.

\mn 2.
\texttt{N.Berestycki@statslab.cam.ac.uk}. University of Cambridge. CMS -- Wilberforce Road, Cambridge, CB3 0WB, United Kingdom.


\newpage

\section{Introduction}

The goal of this paper is to describe the macroscopic behaviour of a process which locally behaves like a Brownian motion, but which is constrained to satisfy a global constraint of a self-avoiding nature. Informally speaking, we consider one-dimensional Brownian motion conditioned on the event $\cE$ that the local time of the process is bounded by a fixed constant, say 1, at every time and position. The event $\cE$ has of course zero probability, so a precise definition is needed - this is deferred to the next section. For the moment, it suffices to say that it is possible to define a probability measure $\mathbb{Q}$ on continuous paths corresponding to this conditioning, which is obtained by a limiting procedure.

\medskip From an intuitive point of view, one expects that, conditionally on $\cE$, the process will be transient and must in fact escape to infinity with a positive velocity. In fact, one expects the speed to be at least equal to 1, since that is precisely what it means to spend less than one unit of local time per level. This being a very costly behaviour for Brownian motion, it is tempting to believe that the process is not likely to satisfy any constraint that would be even stronger, and hence that the speed of the process will in fact be equal to 1 in the limit.

\medskip Our main finding in this paper is that this intuition turns out to be erroneous. To be precise, we obtain the following result.
\begin{theorem} \label{Tballintro}
$
\lim_{t\to \infty}{X_t}/t = \gamma_0
$
exists in $\mathbb{Q}$-probability, and furthermore:
\begin{equation}\label{speed}
\gamma_0 =\frac3{1-2j_0^{-2}}=4.5860 \ldots
\end{equation}
where $j_0=2.4048\ldots$ is the first nonnegative zero of the Bessel function $J_0(x)$ of the first kind and of order 0.
\end{theorem}

A more precise result is stated in Theorem \ref{T:ballistic} in the next section, after the definitions have been given. In Theorem \ref{T:rw}, we obtain a similar result for the corresponding random walk problem: given some $L_0\ge 2$, a simple symmetric random walk on $\mathbb{Z}$ is conditioned to never visit any site more than $L_0$ times. Under the limiting measure $\mathbb{Q}$, we show that the particle escapes to infinity with a certain speed $\gamma(L_0)$ and we show that $\gamma(L_0)>1/L_0$.

\medskip We call this phenomenon Brownian entropic repulsion, by analogy with a situation arising in the study of the harmonic crystal which will be described below (in Section \ref{S:relwork}). Roughly speaking, entropic repulsion describes the fact that the easiest way to achieve a certain global constraint for a random process is to achieve much more than required. Here, this phenomenon arises due to the fact the local time process has wild oscillations, and therefore the process must on average have a small local time if it wants to avoid ever being equal to 1. As discussed in Section \ref{S:relwork}, the situation in the harmonic crystal is not much different. We also describe other conditionings of Brownian motion where a similar entropic repulsion occurs in the paper \cite{bb2}, and recall some results of that paper later on in Section \ref{S:relwork}. We expect entropic repulsion to be a general principle in this sort of situations, even though it seems hard to even formalize this idea precisely.

\medskip Our techniques are very different in the continuous and the discrete case. In the continuous case, our main tools are the Ray-Knight theorem and some careful coupling estimates. The existence and uniqueness of the measure $\mathbb{Q}$ is obtained by showing that the approximating sequence forms a Cauchy sequence for a suitable metric (and showing that this implies weak convergence). The value of the speed is obtained through a connection to an eigenvalue problem for the Laplacian in the unit disk of the plane (Sturm-Liouville problem). While it seems possible to adapt these techniques to the discrete case, we have used here a rather different method which we believe sheds additional light onto the problem. In particular, the notion of regenerating levels plays a significant role in this proof and we crucially apply the renewal theorem in the spirit of Kesten \cite{kesten}. This theorem is here viewed as a purely analytic result on sequences of numbers satisfying certain conditions, and is applied to sequences which do not have obvious probabilistic interpretations.

\section{Statement of the results}


\medskip Let $\Omega$ be the space of continuous, real-valued functions defined on $[0,\infty)$. We endow $\Omega$ with the Skorokhod topology and the Borel $\sigma$-field defined by it, and with the Wiener measure $\bm$. (We let $\bm_x$ be the Wiener measure started at the point $x \in \R$).  Let $(X_t,t\ge 0)$ be the canonical (coordinate) process on $\Omega$, and let $L(t,x)$ denotes a jointly continuous version of the local times process of $X$, i.e., $\bm$-almost surely for all $x \in \mathbb{R}$ and all $t\ge 0$,
\begin{equation}
L(t,x)=\lim_{\eps\to 0} \frac1{2\eps} \int_0^t \indic{|X_s - x | \le \eps}ds.
\end{equation}
(We may, occasionally, write $L(t,x,\omega)$ to make explicit the dependence of $L(t,x)$ upon the path $\omega \in \Omega$). In particular, $X$ satisfies the occupation formula: almost surely, for all $t\ge 0$ and for all nonnegative Borel function $f$,
\begin{equation}\label{occ}
\int_0^t f(X_s)ds = \int_\R f(x) L(t,x)dx.
\end{equation}
(For this and other basic facts about local times of Brownian motion, we refer the reader to Chapter VI of \cite{revuz-yor}. The above statement corresponds to Corollary (1.6) in that reference.) For all $a>0$, let
$$
\tau_a= \inf\{t>0: X_t \ge a \}.
$$
We approximate the event $\cE$ described in the introduction by conditioning on what happens up to time $\tau_a$, and let $a$ tend to infinity. Hence, define
\begin{equation}\label{cE_a}
{\cal E}_a:= \left\{\sup_{t\le \tau_a}\sup_{x\in \R} L(x,t) \le 1\right\}.
\end{equation}

A more precise statement of Theorem \ref{Tballintro} follows. Recall that $\gamma_0 =4.5860\ldots $ is defined by (\ref{speed}).

\begin{theorem} \label{T:ballistic} The family of measures $\{\mathbb{P}_a:=\bm(\cdot|\cE_a)\}_{a\ge 0}$ converges weakly to a measure $\mathbb{Q}$ on $\Omega$ as $a \to \infty$. Moreover,
$
\lim_{t\to \infty}{X_t}/t = \gamma_0
$
in $\mathbb{Q}$-probability. \end{theorem}

Roughly speaking, the idea for the proof of this theorem is that, by the Ray-Knight theorem, the local times of Brownian motion between the two endpoints $0$ and $a$ behave as the square radius of a two-dimensional Brownian motion. Conditioning by the event $\cE_a$ amounts to conditioning this Brownian motion to stay within the unit ball. By well-known results due to Pinsky \cite{pinsky} on \emph{metastability}, this has a simple equilibrium distribution, under which the square radius has an average $c<1$. The asymptotic speed of the process is then given by $\gamma_0 = 1/c$, and thus $\gamma_0 >1$.

\begin{rmk}
\emph{If one requires the local time to be bounded by $C>0$ rather than one in the events $\cE$ and $\cE_a$, it can be shown that the limiting speed of the process becomes $\gamma_0/C$. That is, entropic repulsion makes the particle travel $4.5860\ldots$ as fast as one would expect from the constraint in the conditioning.}
\end{rmk}

\begin{rmk}
  \emph{Conditioning on the event $\cE_a$ drives the process to $+\infty$ because the condition is imposed at time $\tau_a$, the hitting time of $a>0$. If we replace $\tau_a$ in the definition of $\cE_a$ by $\tau'_a$, where
  $$
  \tau'_a = \inf\{t \ge 0:|X_t|=a\}
  $$
  then Theorem \ref{T:ballistic} easily implies that the same statement is true with $\lim_{t \to \infty} X_t /t = \pm \gamma_0$ with probability $1/2$ each.}
\end{rmk}

We now turn to the discrete version of the problem, about which we know both more and less information. As our basic probability space we take $\Om = \{-1, +1\}^{\Bbb
Z_+}$. A generic point of $\Om$ is written as $\om = \{\om_t\}_{t
\ge 0}$. For $\omega \in \Omega$, let
$$
S_n(\omega)= \sum_{j=1}^n \omega_j, \ \ \ \ n = 0,1, \ldots,
$$
be the random walk on $\Bbb Z$ associated to $\omega$.
For $t \in \Bbb
Z_+, x \in \Bbb Z$, define
$$
L(t,x) = \sum_{j=1}^t \indic{S_j=x}.
$$
Of course, $L(t,x)$ is a function on $\Om$. Occasionally it will
be useful to write $L(t,x,\om)$ for the value of $L(t,x)$ at the
point $\om \in \Om$. In fact, $L(t,x,\om)$ depends only on the
first $t+1$ coordinates of $\om$, so we can also regard it as a
function on $\Om_t := \{-1, +1\}^{t+1}$. If $\om_t \in \Om_t$ we
shall also use the notation $L(t,x, \om_t)$ for the value of
$L(t,x)$ at this point. Unless otherwise indicated we take $S_0 =
0$. Let us now define the event $\cB$ which serves as our constraint: for $r \in \mathbb{Z}_+$, let
\begin{align*}
\cB_r &= \{L(\tau_r,x) \le L_0 \text{ for all } x\}\\
&= \{\omega \in \Om: L(\tau_r(\omega),x,\omega) \le L_0 \text{ for all } x\}
\end{align*}
where
$$
\tau_r=\inf\{k\ge 1: S_j =r\}.
$$

\medskip To formulate our main result for random walks we will need to introduce the notion of ``regenerating levels", to borrow from the terminology of random walks in random environments. Define:
$$
\nu_1= \inf\{r\ge 0: S_{t}> r \text{ for all } t\ge \tau_r+1\}
$$
and define recursively, for all $i\ge 2$,
$$
\nu_i= \inf\{r>\nu_{i-1}: S_{t}> r \text{ for all } t\ge \tau_r+1 \}.
$$
The levels $\nu_i$ are those which are visited only once for a given trajectory $\omega$.

\begin{theorem} \label{T:rw} The measures $\mathbb{P}(\cdot | \cB_r)$ converge weakly to a limiting measure $\mathbb{Q}$ as $r\to \infty$.
Then for all $j\ge 1$, $\nu_j <\infty$, $\mathbb{Q}$-a.s. Moreover, the random variables $(\nu_{j+1}-\nu_j)_{j\ge 1}$ are i.i.d. and satisfy
$$
E_\mathbb{Q}(\nu_{j+1}-\nu_j) <\infty.
$$
The portions of the path between two successive renewal levels are also independent. In particular, $\gamma(L_0)=\lim_{k \to \infty} X_k/k $ exists $\mathbb{Q}$-almost surely. and is a nonrandom number satisfying $\gamma(L_0)>1/L_0$.
\end{theorem}

\section{Related work}
\label{S:relwork}

\subsection{Harmonic crystal with hard wall repulsion}

As already mentioned above, the term ``entropic repulsion" was introduced to describe a situation arising in the study of the discrete Gaussian free field on a lattice (also known as the harmonic crystal) with hard wall repulsion, which presents some strong analogy to the phenomenon described by Theorems \ref{T:ballistic} and \ref{T:rw}. Indeed, in \cite{bdg}, the following result (among other things) is proved. Let $\Phi_N =(\phi_x)_{x \in V_N}$ be the law of a free field on a box $V_N=\{1,\ldots,N\}^2$ with zero boundary conditions and covariance $\cov(\phi_x,\phi_y)=G_N(x,y)$ (the discrete Green function stopped when the walk reaches the outside of the box). Let $D_N$ be a ``nice" domain in the box (essentially, the discrete approximation of a smooth fixed domain in $(0,1)^2$ away from the boundary, blown up by a factor of $N$), and let $\Omega^+_{D_N}$ be the event that $\phi_x \ge0$ for all $x \in D_N$. Then, conditionally on $\Omega^+_{D_N}$, the value of the field $\phi_x$ is typically of order $\log N$, in the following strong sense: for all $\eps>0$,
\begin{equation}
\lim_{N\to \infty} \sup_{x \in D_N} \mathbb{P}\left(\left.|\phi_x - \frac4\pi \log N|\ge \eps \log N \right|\Omega^+_{D_N}\right)=0
\end{equation}
The intuitive reason for this behaviour is the same as in Theorem \ref{T:ballistic} above. To simplify, due to the wild oscillations of the free field (or the local time field, in our case), the simplest way to achieve the constraint is a global shift which guarantees that the wild oscillations do not break the constraint.

\subsection{Brownian motion with limited local time}

\bigskip In \cite{bb2}, we have also studied other conditionings of Brownian motion which favor a self-avoiding behaviour, even though the constraint is much softer than the event $\cE$. Namely, we discuss Brownian motion conditioned on the event $\mathcal{K}$ that the growth of the local time at the origin, is slower than some function of time $f(t)$ where $f$ is nondecreasing but $f(t)t^{-1/2}$ is nonincreasing. We show that if $\int_{1}^\infty f(t)t^{-3/2}dt < \infty$ then the process is transient. We believe this condition to be sharp. In particular, if $f(t) \sim c \sqrt{t} (\log )^{-\gamma}$ for some $c>0$ and $\gamma \ge 0$, then the process is transient as soon as $\gamma >1$. In the regime where $0\le \gamma \le 1$, and where we thus anticipate that the process is recurrent, we nonetheless expect an entropic repulsion phenomenon to occur in the sense that $L_t = o(f(t))$ with high probability for $t \to \infty$.

\subsection{Edwards and Domb-Joyce polymer models}

Finally, the present work is closely related to the field of polymer models. The well-known Domb-Joyce model (and its Brownian analogue -- the Edwards model) is a model where simple random walk measure is penalized by a weight exponential in the number of self-intersections. More precisely, given an inverse temperature $\beta>0$, the Domb-Joyce model is defined by looking at the measure $\mu_N$ on nearest-neighbours discrete random paths of length $N$ obtained by setting
$$
\mu_N(\omega)=\frac1{Z_N}\exp\left(-\beta\sum_{0\le i <j \le N}\indic{\omega_i=\omega_j}\right) 2^{-N}
$$
where $Z_N$ is a normalizing constant. Similarly, the Edwards model (in one dimension) is defined by taking a large $T>0$ and considering the measure $\mu_T$ whose density with respect to the Wiener measure is
\begin{equation}\label{westwater}
\frac{d\mu_T}{d\mathbb{W}}= \frac1{Z_T}\exp\left(-\beta \int_{\mathbb{R}}L(T,x)^2dx\right).
\end{equation}
where $L(t,x)$ is a jointly continuous version of the local time at time $t$ and position $x$.
It is the limit of the distribution of the position of the endpoint under these measures (and their dependence on $\beta$) as $N$ or $T$ tend to infinity which is of interest. The main result on this model, proved in \cite{hhk}, is that $X_T$ is approximately normally distributed with a mean $c(\beta)T$ and variance $\sigma^2(\beta)T$. In the case of the Edwards model, these parameters have simple dependencies on $\beta$: in fact, the variance parameter $\sigma^2(\beta)$ is independent of the self-repellency strength $\beta$, while $c(\beta)=b^*\beta^{1/3}$ for some $0<b^*<\infty$. However, in the discrete Domb-Joyce model, the dependency on $\beta$ is largely unknown -- it is still an open question to show that $c(\beta)$ is monotone in $\beta$. See \cite{polreview} and the references therein for a very interesting account of the theory. See also the paper \cite{nrw} for a polymer model related to our work, where explicit calculations on the ballistic behaviour of the process can be done.

We note that both the present work (in the continuous case) and the papers \cite{hhk,nrw} use in a fundamental way the Ray-Knight theorem, as well as (for \cite{hhk}) a connection to an eigenvalue problem for the Laplacian. However, this is where the analogy stops: while \cite{hhk} requires many difficult analytical estimates, we only require careful but simple-minded probabilistic coupling estimates. Also, in this paper we discuss the full convergence of the path $(X_s, s \ge 0)$ (in the sense of weak limits of measure on paths) rather than its position at a large time. During the revision of this paper we learnt from the referee that Joseph Najnudel \cite{najnudel} has recently constructed a probability measure on $\Omega$ corresponding for the whole process in the setup of the Domb-Joyce polymer model. His techniques are very different from ours. Note also that the discrete case uses entirely different techniques. Finally, we mention that it is very likely that our techniques would yield a central limit theorem for the position of the particle in Theorem \ref{T:ballistic} and \ref{T:rw}. We have not tried to pursue this direction.

A related problem has also been studied by M\"orters and Sidorova \cite{MortersSidorova}, where they analyse the order of magnitude of the maximal displacement of a random walk conditioned on the $p^\text{th}$ moment of its local time profile being unusually small, for some $p>1$. More precisely, let
$$
\Lambda_n(p) = \sum_{z \in \mathbb{Z}} L(z,n)^p,
$$
where $L(z,n)$ denotes the the number of visits by a simple random walk to $z$ by time $n$. They consider the simple random walk conditioned on the event that $\{\Lambda_n(p) < \eps_n E(\Lambda_n(p))\}$ for some sequence $\eps_n = o(1)$. They are able to show that under this conditioning, there exist constants $c_1, c_2>0$ such that
$$
c_1 \le \frac{\max_{1 \le i \le n} |S_n|}{ \sqrt{n} \eps_n^{-1/(p-1)}} \le c_2
$$
with high probability as $n\to \infty$. Their result is based on a careful large deviation analysis.

\section*{Acknowledgements}

All our gratitude goes to Harry Kesten, whose help for Theorem \ref{T:rw} was invaluable.
We thank A.S. Sznitman for pointing out the connection to the work on the entropic repulsion of the harmonic crystal, and Ross Pinsky for pointing out several typos in a first version of this paper. We also thank an anonymous referee for several comments which improved the paper.

\section{Existence and uniqueness of the weak limit.}

We start the proof of Theorem \ref{T:ballistic} with the existence and uniqueness of a weak limit for the measures $\mathbb{P}_t:= \bm( \cdot | \cE_t)$ as $t \to \infty$, for the Skorokhod topology (we refer the reader to \cite{billingsley} for background on weak convergence). In fact, we are going to prove a stronger statement and show that for the total variation distance on sets measurable with respect to $\mathcal{F}_{\tau_a}$ for a fixed arbitrary $a>0$, the measures $\mathbb{P}_t$ form a Cauchy sequence. For the convenience of the reader, we first explain precisely what we mean by this and then prove that this implies weak convergence with respect to the Skorokhod topology. The remainder of the section will be devoted to the proof that $\mathbb{P}_t$ is a Cauchy sequence in that sense.

Thus, let $a>0$, and recall that $\tau_a= \inf\{s>0: X_s \ge a\}$. For probability measures $\mu, \nu$ on $(\Omega, \cF_{\tau_a})$ define:
\begin{equation}\label{TV}
d_a( \mu ,\nu ) := \sup_{A \in \cF_{\tau_a}} | \mu(A)-\nu(A)|.
\end{equation}
\begin{lemma}\label{L:cauchy}
Let $\{\mu_t\}_{t\ge 0}$ be a sequence of probability measures on $\cF$ such that for every $a$, the restrictions of $\mu_t$ to $\cF_{\tau_a}$ forms a Cauchy sequence for the distance $d_a$, i.e., for every $\eps>0$, there exists $t_0$ such that for all $s,t\ge t_0$,
\begin{equation} \label{cauchy}
d_a(\mu_t, \mu_s) \le \eps.
\end{equation}
Assume also that for all $A>0$ fixed, \begin{equation}\label{22}
\lim_{b \to \infty} \limsup_{t \to \infty} \mu_t\left(\sup_{s < A}
|X(s)| > b\right) = 0.
\end{equation}
Then there exists $\mu$ a probability measure on $(\Omega, \cF)$ such that $\mu_t \to \mu$ weakly as $t \to \infty$ for the Skorokhod topology on $\Omega$.
\end{lemma}

\begin{proof}
The proof is mostly routine manipulations, so we content ourselves with outlining it. The bottom line is that convergence in total variation distance is typically much stronger than weak convergence. Fix $A \in \cF_{\tau_a}$. Then by (\ref{cauchy}) and (\ref{TV}), we get that $\mu_t(A)$ is a Cauchy sequence, so has a limit $\mu^a(A)$ as $t \to \infty$. It is easy to check that $\mu^a(A)$ is a probability measure on $(\Omega, \cF)$ (the $\sigma$-additivity property follows from the uniformity over all sets in (\ref{TV}). Moreover, we get that for every $A \in \cF_{\tau_a}$, $\mu_t(A) \to \mu^a(A)$ as $t \to \infty$. From this it is trivial to check that $\mu^a$ satisfies the conditions of Kolmogorov's extension theorem, and thus we may define a unique measure $\mu$ such that for all $a>0 $,
\begin{equation}
\label{conv}
 \mu_t(A) \underset{t \to \infty}{\longrightarrow} \mu(A), \ \  \text{ for all } A \in \cF_{\tau_a }.
\end{equation}
While it does not seem \emph{a priori} easy to extend (\ref{conv}) to all sets $A \in \cF$ such that $\mu(\partial A)=0$, where $\partial A$ is the boundary of the set $A$ with respect to the Skorokhod topology, we claim that it follows easily from (\ref{cauchy}) that $\{\mu_t\}_{t\ge 0}$ is a tight family. There are two conditions to verify, of which the first one (non-explosion in finite time) is part of the assumption on $\mu_t$ (see (\ref{22})). the second condition to verify is: for all $A>0$, and for each $\eta > 0$
\begin{equation}\label{23}
\lim_{\ep \downarrow 0} \limsup_{t \to \infty}
\mu_t\big(\sup\{|X(s')-X(s'')|: 0\le  s',s'' < A, |s'-s''|< \ep\} >
\eta\big) = 0. 
\end{equation}
First observe that by (\ref{22}) and by (\ref{conv}), we also have that for all $s>0$, and for all $E \in \cF_s$, \begin{equation}\label{conv2}
\mu_t(E) \to \mu(E)
\end{equation}
uniformly in $E \in \cF_t$ as $t \to \infty$. Therefore, fix $\delta>0$ and let $t_0$ be such that for all $E \in \cF_A$, $|\mu_t(E) - \mu_s(E)| \le \delta$ for all $t,s \ge t_0$. Thus for $t \ge t_0$, and $E_{\eta, \ep}$ the event in (\ref{23})
$$|\mu_t(E_{\eta, \ep}) - \mu_{t_0}(E_{\eta, \ep}) | \le \delta$$
from which it follows that
$$
\limsup_{t \to \infty} \mu_t(E_{\eta, \ep}) \le \delta + \mu_{t_0}(E_{\eta, \ep}).
$$
Now, since $\mu_{t_0}$ is the law of a continuous process, $\lim_{\ep \to 0} \mu_{t_0} (E_{\eta, \ep}) = 0$. Therefore,
$$
\lim_{\ep \downarrow 0} \limsup_{t \to \infty}
\mu_t (E_{\eta, \ep}) \le \delta
$$
where $\delta>0$ is arbitrary. (\ref{23}) follows by letting $\delta \to 0$. Therefore, $\{\mu_t\}_{t \ge 0}$ forms a tight family, and so there exists some weak subsequential limit. On the other hand, by (\ref{conv2}) this limit must be $\mu$ since the finite-dimensional marginal distributions are specified by events of the form $E \in \cF_t$ for some finite $t>0$. Since the weak subsequential limit is unique and we have proved tightness, we conclude that $\mu_t \to \mu$ weakly as $t \to \infty$, for the Skorokhod topology.
\qed
\end{proof}

\medskip For the proof that there exists a weak limit to the sequence $\bm(\cdot | \cE_a)$ as $a \to \infty$, we will use Lemma \ref{L:cauchy}. It turns out that (\ref{22}) is very easy to verify, and the core of the proof is to check (\ref{cauchy}). Crucial to this proof is the Ray-Knight theorem; we start by reminding the reader the statement of this result, as can be found in \cite{revuz-yor}, chapter XI.2, or \cite{rw} (VI. (52.1)) for the formulation we use here.

A square Bessel process of dimension $\delta \ge 0$ is the unique strong solution to the stochastic differential equation:
\begin{equation}
Z_t = z_0 + 2\int_0^t \sqrt{|Z_s|} dB_s + \delta t, \ \ \ z_0\ge 0 .
\end{equation}
In the special case where $\delta=0$ this process is known as the Feller diffusion. When $\delta$ is an integer $\ge 1$, $Z$ can be interpreted as the square Euclidean norm of a $\delta$-dimensional Brownian motion.

Let $(B_t,t\ge 0)$ be a one-dimensional Brownian motion with joint local time process $\{L(t,x)\}_{t\ge 0, x \in \R}$, and let $\tau_a=\inf\{t\ge 0: B_t =a\}$ be the hitting time of a fixed level $a>0$.

\begin{theorem} \label{T:rk} \text{\em (Ray, Knight)} For all $a>0$, the law of $L(\tau_a, a-x)$ is specified by:

\begin{enumerate}

\item $\{L(\tau_a,a-x)\}_{0\le x \le a}$ is a square Bessel process of dimension 2, started at 0.

\item Conditionally given $L(\tau_a,0)=z_0\ge 0$, $\{L(\tau_a,-x)\}_{x\ge 0}$ is a Feller diffusion started at $z_0$ and is independent from $\{L(\tau_a,a-x)\}_{0\le x \le a}$.

\end{enumerate}

\end{theorem}

We now state a lemma which allows us to compare different constraints on 2-dimensional square Bessel process, which will be used repeatedly throughout the proof. It should be noted that in general there is no known way to compare the effect of two different constraints, even when one is intuitively stronger than the other. Lemma \ref{Yf-ub} shows however that making comparisons is possible when, in some sense, the constraint only deals with the position of the process.

 Let $(Y_t,t\ge 0)$ denote a square Bessel process of dimension 2. We may view $Y$ as a random element of $(\Omega, \cF)$ under the probability measure $\be$, which is (as explained above) the law on $(\Omega, \cF)$ of the squared Euclidean norm of a two-dimensional Brownian motion. As we work with many different processes it will at times be convenient to use a generic symbol $P$ for the underlying probability space of different random processes. The notations $\{X(t)\}_{t\ge 0}$, $\{Y(t)\}_{t\ge 0}$ then serve to differentiate these processes, and from the context it should be clear to which processes they refer.

For a given $T>0$ and a positive measurable function $f:[0,T] \to [0,\infty)$, let $\{Y^f(t)\}_{t\ge 0}$ denote a version of $Y$ conditionally given $\cA(0,x)=\{\omega \in \Omega: \omega_s \le f(s) \text{ for all }s\le T; \text{ and } \omega_0 =x\}.$ For $0\le u \le T$, we also define the event
\begin{equation}\label{A(u,x)}
\cA(u,x)=\{\omega \in \Omega: \omega_s\le f(u+s) \text{ for all } s\le T-u ; \text{ and } \omega_0=x\}.
\end{equation}

\begin{lemma} \label{Yf-ub}
$\{Y^f(t),0\le t\le T\}$ is an inhomogeneous diffusion on $\mathbb{R}_+$ which satisfies:
\begin{equation} \label{Yf}
dY^f(t)=\sqrt{Y^f(t)}dB_t + \{2-\delta^f(t,Y^f(t))\}dt
\end{equation}
where $\delta^f(t,y)\ge 0$ for all $0\le t\le T$ and for all $0<y<f(t)$. Moreover, if $g$ is another function such that $g(t) \le f(t)$ for all $0\le t \le T$, then
\begin{equation} \label{Yg}
\delta^g(t,y) \ge \delta^f(t,y)
\end{equation}
for all $0<y <g(t)$. As a result, $Y \succeq Y^f \succeq Y^g$, where $\succeq$ stands for stochastic domination.
\end{lemma}

\begin{proof}
It is a well-known fact the conditioned process $Y^f$ can be realized as an $h$-transform of the original process $Y$: more precisely, by Girsanov's theorem, $Y^f$ is an inhomogeneous diffusion having the form (\ref{Yf}) where
\begin{equation} \label{delta}
\delta^f(t,y) = -\frac{\partial}{\partial y} \log h(t,y)
\end{equation}
where
\begin{equation}\label{h}
h(t,y)= \be(\cA(t,y)).
\end{equation}
(Details can be found for instance in \cite{rw}, IV.39 in the case where the process $Y$ is Brownian motion. Generalization to weak solutions of stochastic differential equations presents no difficulty and we do not give the details here).
For the first part of Lemma \ref{Yf-ub} it thus suffices to prove that
\begin{equation}
\frac{\partial h(t,y)}{\partial y}\le 0.
\end{equation}
Let $\eps>0$, and let $y=x+\eps$. It suffices to prove that for all $\eps>0$ small enough,
\begin{equation} \label{AxAy}
\be (\cA(t,x)) \ge \be(\cA(t,y)).
\end{equation}
We use a coupling technique to prove this. Let $Y_1$ denote a square Bessel process of dimension 2 started from $x$ and let $Y_2$ denote an independent square Bessel process of dimension 2, but started from $y$. Let $t>0$ and let
\begin{equation}\label{fhat}
\hat f (s) = f(T-t+s)\text{  for all }s\le \hat T:= T-t.
\end{equation}
Let $\tau=\inf\{t>0: Y_1(t)=Y_2(t)\}$, and let
\begin{equation}
Y_3(t)=\begin{cases}
Y_1(t) &\text{ if }  t< \tau\\
Y_2(t) & \text{ else.}
\end{cases}
\end{equation}
Then by the strong Markov property, $Y_3$ has the same distribution as $Y_1$ and moreover $Y_3(s) \le Y_2(s)$ for all $s\ge 0$ almost surely. It follows that if
$$
Y_2(s) \le \hat f(s) \text{ for all $s \le \hat T$}
$$
then automatically
$$
Y_3(s)\le \hat f(s) \text{ for all $s \le \hat T$}.
$$
The desired (\ref{AxAy}) follows.

The second part of Lemma \ref{Yf-ub} is an easy consequence of the first part. Indeed, $Y^g$ can be obtained by conditioning further the process $Y^f$ to stay below the function $g$. We conclude again by Girsanov's theorem that there exists an additional drift term $\delta^{f,g}(t,y)$ such that
\begin{equation}
dY^g(t)= \sqrt{Y^g(t)}dB_t + \{2 -\delta^f(t,Y^g(t)) - \delta^{f,g}(t,Y^g(t))\}dt
\end{equation}
and that $\delta^{f,g}(t,y)$ satisfies:
\begin{equation}
\delta^{f,g}(t,y)= \frac{\partial}{\partial y} \log h^{f,g}(t,y).
\end{equation}
This time,
\begin{equation}
h^{f,g}(t,y) = P(Y^{\hat f}(s) \in \cA'(t,y))
\end{equation}
where $\hat f$ is defined in (\ref{fhat}) and $\cA'(t,y)$ has the same definition as $\cA(t,y)$ except $f$ is replaced with $g$. Since $Y^{\hat f}$ is a strong Markov process by the first part, the coupling argument works equally well to show that
\begin{equation}
\frac{\partial}{\partial y} h^{f,g}(t,y) \le 0.
\end{equation}
As above, this implies $\delta^{f,g}(t,y)\ge 0$ for all $t\le T, 0<y < g(t)$, and thus $\delta^g (t,y)\ge \delta^f(t,y)$. To get the final statement of the lemma, we note that it is easy to show that (\ref{Yf}) admits strong and pathwise unique solutions, since the coefficients are locally Lipschitz. From this and Theorem 3.7 in \cite{revuz-yor}, the desired stochastic dominations follow directly.
\qed
\end{proof}

\medskip We now show that if $\mathbb{P}_t:= \bm( \cdot | \cE_t)$ for $t>0$, then $\mathbb{P}_t$ satisfies the assumptions of Lemma \ref{L:cauchy}. Let $s,t>0$ with $s<t$ and let $0< a < s$. First note that if $A \in \cF_{\tau_a}$, then we have by elementary manipulations:
\begin{equation}\label{coupl0}
\mathbb{P}_t(A) = \mathbb{P}_s(A) \frac{\bm( \cE_t | A \cap \cE_s)}{\bm( \cE_t | \cE_s)}
\end{equation}
so it suffices to prove that the ratio is arbitrarily close to 1, uniformly in $s,t$ large enough, and $A \in \cF_{\tau_a}$. We will show the existence of a coupling $P$ between two processes $X$ and $Y$, having respectively the law of $\bm(\cdot | \cE_s)$ and $\bm(\cdot | \cE_s \cap A)$ such that $P$-almost surely,
\begin{equation}\label{coupl1}
L(\tau_s, x, X)= L(\tau_s,x,Y) \text{ for all $x \ge a+\Delta$}
\end{equation}
where $\Delta <\infty$ $P$-almost surely, and in fact there exists $\Delta^*$ a random variable whose distribution does not depend on any parameter, and such that $\Delta \preceq \Delta^*$, and $\Delta^*< \infty$ almost surely. For the moment, let us admit these facts and see how we proceed with them. Let $(Z_u, u \ge 0)$ be a square Bessel-0 process started at an unspecified point $Z_0=x \in (0,1)$. We claim that there exists $C,\alpha >0$ independent of $x$ such that for all $t>0$:
\begin{equation}
\label{exp0}
P\left(Z_t >0 \left| \sup_{0\le u } Z_u \le 1 ; Z_0 =x\right. \right) \le C e^{- \alpha t}.
\end{equation}
This follows easily from the Markov property and the fact in any period of duration 1, $Z$ started from position 1 has a positive probability $p_0$ to reach zero. If not, then at the next iteration the process is still below 1 and again has a probability bigger than $p_0$ to die out in the next interval. Using Lemma \ref{Yf-ub}, we conclude that (\ref{exp0}) holds with $\alpha = - \log(1-p_0)$.

To ease notations, let $F_1(u) = 1- L(\tau_s, t-u,X)$ for all $u\ge 0$, and similarly let $F_2(u) = 1- L(\tau_s, t-u,Y)$ for all $u\ge 0$. In particular, note that $F_1(u)=F_2(u)=1$ for all $u \le t-s$, $P$-almost surely. Note that our assumption (\ref{coupl1}) implies that $F_1(u) = F_2(u) $ for all $u \le t-(a+\Delta)$. If now $(Z_u,u\ge 0)$ is the Ray-Knight diffusion changing dimension at time $u=t-s$, then we have:
$$
\bm( \cE_t | \cE_s) = P( Z_s \le F_1(s) \text{ for all }s \ge 0)
$$
and
$$
\bm( \cE_t | A \cap \cE_s) = P( Z_s \le F_2(s) \text{ for all }s \ge 0).
$$
Let $E_1, E_2$ be the two events in the above equations. It follows that if $p:=\bm( \cE_t | \cE_s)=P(E_1) $ and $q:=\bm( \cE_t | A \cap \cE_s)=P(E_2)$, we have:
\begin{align*}
p&= P(E_1; Z_{t-a-\Delta}=0) + P(E_1; Z_{t-a-\Delta}>0)
\end{align*}
while:
\begin{align*}
q&= P(E_2; Z_{t-a-\Delta}=0) + P(E_2; Z_{t-a-\Delta}>0).
\end{align*}
By definition of $\Delta$, we must have $P(E_1; Z_{t-a-\Delta}=0) = P(E_2; Z_{t-a-\Delta}=0) $, so it follows:
\begin{align*}
|p-q|&= P(E_1; Z_{t-a-\Delta}>0) + P(E_2; Z_{t-a-\Delta}>0).
\end{align*}
and thus
\begin{equation}\label{coupl2}
\left|1-\frac{q}p\right|= P(Z_{t-a-\Delta}>0|E_1) + \frac{P(E_2; Z_{t-a-\Delta}>0)}p.
\end{equation}
We study the two terms in the right-hand side separately. For the first term, we note that by (\ref{exp0}) and Lemma \ref{Yf-ub}, we get
\begin{align*}
P(Z_{t-a-\Delta}>0|E_1) & \le E( Ce^{-\alpha (s-a-\Delta)_+}) \\
&  \le C \  E( e^{-\alpha (s-a-\Delta^*)_+})
\end{align*}
where $x_+=\sup(x,0)$ is the positive part of $x$. Similarly, the second term in (\ref{coupl2}) satisfies:
\begin{align*}
\frac{P(E_2; Z_{t-a-\Delta}>0)}p& = \frac{P(E_2; Z_{t-a-\Delta}>0)}q \frac{q}p \\
&= P(Z_{t-a-\Delta}>0|E_2) \frac{q}p \\
&  \le C \  E( e^{-\alpha (s-a-\Delta^*)_+}) \frac{q}p
\end{align*}
by another application of Lemma \ref{Yf-ub}. To put these two things together, define $\eps:= C \  E( e^{-\alpha (s-a-\Delta^*)_+})<\infty$ and let $x= q/p$. Thus we have proved:
$$
|1-x| \le \eps + \eps x.
$$
Thus $x-1 \le \eps + \eps x$ and solving this inequality we find $x \le (1+\eps)/(1-\eps) = 1 +2 \eps/(1-\eps)$.
Note that by the Lebesgue convergence theorem, as $s \to \infty$, $\eps \to 0$. Thus if $s$ is large enough that $\eps / (1- \eps) \le 2 \eps$, we have proved that
$$
x \le 1 + 4 \eps
$$
and a similar lower bound follows without any difficulty. From this and (\ref{coupl0}), we obtain that for any $\eta>0$, there exists $s_0>a$ large enough that for all $s_0 < s <t$, we have
$$
|\mathbb{P}_t(A) - \mathbb{P}_s(A)| \le \eta
$$
for all event $A \in \cF_{\tau_a}$. In other words, we have proved that if $\mu_t$ is the law of $(X_r, r \le \tau_a)$ under $\mathbb{P}_t$,
$$
d_a (\mu_t , \mu_s)\le \eta
$$
for all $s,t \ge s_0$. That is, $\{\mu_t\}_{t\ge 0}$ forms a Cauchy sequence for the total variation distance. Condition (\ref{22}) is a direct consequence of Lemmas \ref{L:hit2} and \ref{L:tb}, so the proof is deferred to the next section. Thus, provided (\ref{coupl1}) holds, $\mu_t$ satisfies the assumptions of Lemma \ref{L:cauchy} and therefore has a (unique) weak limit $\mu$.

We now turn to the proof of (\ref{coupl1}). This is based on a time-reversal argument and coupling. Note first that it suffices to construct a coupling of $\{L(\tau_s,x, X)\}_{x \in \R}$ and $\{L(\tau_s,x, Y)\}_{x \in \R}$ which achieves  (\ref{coupl1}). Combining the Markov property at time $\tau_a$ for $X$ with the Ray-Knight theorem, we get that if $Z_x = L(\tau_s,x, X)$, then conditionally  on $ \{Z_a= z \in (0,1)\}$, the process $\{Z_{s-x}\}_{ 0 \le x \ge s-a}$ is a square Bessel-2 process conditioned to never exceed 1 and conditioned to be at $z$ at time $s$. In other words, it is a square Bessel bridge of dimension 2 from 0 to $z$ of duration $s-a$, conditioned never to exceed 1 on that interval. There is naturally a similar description for $Y$: if $Z'_x =\{L(\tau_s,x,Y)\}$, then conditionally on $\{Z'_a = z' \in (0,1)\}$, the process $(Z'_{s-x}, 0 \le x \ge s-a)$ is a square Bessel bridge from 0 to $z'$ in duration $s-a$, conditioned never to exceed 1 during that interval. We can now return time and say that, conditionally on $\{Z_a= z \in (0,1)\}$ (resp. $\{Z'_a=z'\}$), the process $(Z_x, a \le x \ge s)$ (resp. $(Z'_x, a \le x \le s)$) is a square Bessel bridge from $z$ (resp. $z'$) to $0$ in duration $s-a$, conditioned never to exceed 1 during that interval. This being a (inhomogeneous) Markov processes, we can couple the processes $Z$ and $Z'$ after the first time (above level $a$) that they meet. That is, let $z<z' \in (0,1)$ without loss of generality, and let $(Z_x, s\le x \le a)$ and $(Z'_x, s\le x \le a)$ be two independent square bessel bridges conditioned never to exceed 1, started respectively from $z$ and $z'$. Consider $\Delta = \inf\{x \ge a: Z_x = Z'_x\}$. Then the process $\hat Z$ defined by
$$
\hat Z_x = \begin{cases}Z'_x & \text{ if $x \le \Delta$}\\
Z_x & \text{ if $x > \Delta$}
\end{cases}
$$
has the same distribution as $Z'$ and satisfies (\ref{coupl1}). It thus suffices to show that there exists $\Delta^*$ independent of $z,z'$, and $a,s$ and $t$, such that
\begin{equation}\label{dom0}
\Delta \preceq \Delta^*
\end{equation}
and $\Delta^* < \infty$ almost surely. We will show that in any interval of duration 1, the two processes have a positive probability $p$ to meet, independently of anything in their past. This will show the inequality (\ref{dom0}) holds with $\Delta^*$ a certain geometric random variable. By lemma \ref{Yf-ub}, $Z'$ is stochastically dominated by an unconditional square Bessel-2 process started from 1, so for any $s \le x \le a$, and for any $\eta>0$,
$$
P( \inf_{y \in [x, x+1]} Z'_y \le \eta | \sigma(Z'_y,y \le x)) \ge p_1 \ \ , a.s.
$$
for some $p_1>0$ (note that $p_1=p_1(\eta)$ depends only on $\eta$). This provides an upper bound for $Z'$ and it remains to give a similar lower bound for $Z$. This takes a few more steps: indeed, it is not hard to see that by the second part of Lemma \ref{Yf-ub}, for any $x \in [s,a]$, $(Z_y, x \le y \le x+1)$ dominates stochastically a square-Bessel bridge $(b_y, x \le y \le x+ 1)$ of dimension 2 from 0 to 0 in duration 1, conditioned on the event $E=\{ \sup_{x \le y \le x+1} b_y \le 1\}$. This event $E$ has positive probability, $p_2$ say. It follows that
\begin{align*}
P( \sup_{y \in [x, x+1]} Z_y \le \eta | \sigma(Z'_y,y \le x)) & \le P( \sup_{x \le y \le x+1} b_y \le \eta |E) \ \, a.s.\\
&\le P( b_{x+1/2} \le \eta)/p_2.
\end{align*}
Now, as $\eta\to 0$, the right-hand side tends to 0, so we can find $\eta>0$ small enough (and universal) such that the right-hand side is smaller than $1/2$ say. Taking the corresponding $p_1(\eta)$, it follows from the above considerations that
$$
P( \Delta \le x+1 | \Delta \ge  x) \ge p_1(\eta)/2
$$
so taking $\Delta^*$ a geometric random variable with success probability $p_1(\eta)/2$ gives us what we were looking for. \qed

\section{Ballistic behaviour}

We start with the identification of the value of the limiting speed, which is obtained by solving a certain eigenvalue problem for the Laplacian in two dimensions. Here again our main tools are the Ray-Knight theorem and some careful comparisons obtained through coupling arguments.

Let $\be^T$ be the law of a square Bessel process of dimension 2 $(Y_t,t\ge 0 )$, conditioned on $\{\sup_{s\le T} Y_s \le 1\}$. The expectation under this probability measure will be denoted by $E_{\be^T}(X)$ for a random variable $X\ge 0$.
\begin{lemma} \label{L:speed} We have:
\begin{equation}
\lim_{t\to \infty}\lim_{T\to \infty} E_{\be^T}(Y_t)= m_0=\gamma^{-1}_0 =\frac{1-2j_0^{-2}}3
\end{equation}
\end{lemma}

\begin{proof}
\mn Step 1. We start by observing that the measure $\be^T_0$ is the law of $(|Z^T(t)|^2,t\ge 0)$, where $Z^T$ is a 2-dimensional Brownian motion conditioned not to exit the unit disc $\mathbf{D}$ by time $T$. By a theorem of Pinsky \cite{pinsky}, the distribution of $\{Z(t)\}_{t\ge 0}$ converges as $T\to \infty$ to a diffusion $\{Z^\infty(t)\}_{t\ge 0}$, which can be determined explicitly. We will not be interested in the precise form of the generator of $Z^\infty$. However we will need to focus on the long term behaviour of the process $Z^\infty$. From the same paper, it is known that $Z^\infty$ admits an invariant nontrivial probability measure measure $\pi$ on $\mathbf{D}$ whose density is equal to:
\begin{equation}
\pi(dx)= \frac1C\varphi(x)^2 dx
\end{equation}
where 
$\varphi$ is the principal eigenfunction associated with the smallest eigenvalue of the operator $L=-\frac12 \Delta$ with Dirichlet boundary conditions on $\mathbf{D}$, and $C= \int_{\mathbf{D}} \varphi(x)^2 dx $. (Note that it does not matter how we have normalised $\varphi$ here). That is,
\begin{equation}\label{dirichlet}
\begin{cases}
\frac12 \Delta \varphi &= -\lambda \varphi \\
\varphi|_{\partial {\mathbf{D}}}&=0.
\end{cases}
\end{equation}
It is well-known that the problem (\ref{dirichlet}) has solutions only for a discrete set of values $\{\lambda_0 < \lambda_1 <\ldots\}$ where the lowest eigenvalue is simple: i.e., the corresponding eigenspace is one-dimensional, generated by an eigenfunction denoted by $\varphi_0$, the principal eigenfunction. Thus $\varphi=\varphi_0$, which is well-known to be rotationally invariant (a good reference at this level of generality is Jost \cite{jost}, Chapter 9.5). Hence $\varphi(x)$ takes the same value over the entire circle of radius $0<r<1$. We may thus define a function $\phi(r)$ on $(0,1)$ such that $\phi(r) = \varphi(x)$ for all $x \in \mathbf{D}$ such that $|x|=r$. By the ergodic theorem (\cite{rw}, V.54) applied to the diffusion $(Z^\infty_t,t\ge 0)$, it follows that
\begin{equation}
\lim_{t\to \infty} E(|Z^\infty_t|^2)=m_0:= \frac1C\int_{\mathbf{D}} |x|^2\varphi(x)^2dx
\end{equation}
Therefore, $\lim_{t\to \infty}\lim_{T\to \infty}E_{\be^T_0}(Y_t)$ exists and is equal to $m_0=\frac1C\int_{\mathbf{D}} |x|^2\varphi(x)^2dx$.

\mn Step 2. It turns out that this integral can be evaluated explicitly. The principal eigenfunction can be identified explicitly as (see, e.g., Courant and Hilbert \cite{ch} (29) in Chapter V)
\begin{equation} \label{bessel}
\phi(r) =  J_0(j_0 r)
\end{equation}
where $J_0$ is the Bessel function of the first kind of order $\nu=0$, $j_0=2.4048\ldots$ is the first nonnegative zero of $J_0$. Having chosen this normalisation of $\phi$, $C$ is given by
\begin{equation}
C=\int_0^1 J_0(j_0 r)^2 2 \pi r dr.
\label{bess const}
\end{equation}
It follows that
\begin{equation}
\label{m eval0}
m_0=\frac{2\pi\int_0^1 r^3 J_0(j_0 r)^2dr}{2\pi\int_0^1 rJ_0(j_0 r)^2dr} =j_0^{-2} \frac{\int_0^{j_0} x^3 J_0(x)^2dx }{\int_0^{j_0} xJ_0(x)^2dx }.
\end{equation}
We turn to the following result which can be found in \cite[p. 137]{watson}, known as Schafheitlin's reduction formula: for all $z\ge 0$, and all $\mu \ge 0$,
\begin{align}
(\mu+2) \int_0^z \hspace{-0.25cm} x^{\mu+2}J_0(x)^2 dx &= -(1/4)(\mu+1)^3\int_0^z x^\mu J_0(x)^2dx \nonumber \\
&  + \frac12\left[x^{\mu+1}(xJ_0'(x)-\frac12(\mu+1)J_0(x))^2\right.\nonumber \\
& \hspace{.8cm} \left. +x^{\mu+1}(x^2+\frac14(\mu+1)^2)J_0(x)^2\right]_0^z \label{sch}
\end{align}
Taking $\mu=1$ and $z=j_0$ and recalling that $J_0(j_0)=0$, we obtain:
\begin{equation}\label{meval 1}
3\int_0^{j_0} x^3 J_0(x)^2dx = -2\int_0^{j_0} xJ_0(x)^2dx +\frac12j_0^4J_0'(j_0)^2
\end{equation}
Thus
\begin{equation}\label{meval 2}
3\frac{\int_0^{j_0} x^3 J_0(x)^2dx }{\int_0^{j_0} xJ_0(x)^2dx }=-2+ \frac{j_0^4J_0'(j_0)^2}{2\int_0^{j_0} xJ_0(x)^2dx}.
\end{equation}
It also turns out that
\begin{equation} \label{meval 3}
\int_0^{j_0} xJ_0(x)^2dx =\frac12j_0^2 J'_0(j_0)^2.
\end{equation}
(This is a consequence of the fact that the Bessel functions are orthonormal for the weight $x$: this is a classical property which can be found in \cite[p. 576]{watson} for instance). 
Thus, 
using (\ref{m eval0}) together with (\ref{meval 2}) and (\ref{meval 3}) we obtain:
\begin{equation}
\label{meval 5}
m_0=j_0^{-2}(1/3)[-2+j_0^2] = (1/3)[1-2j_0^{-2} ].
\end{equation}
This completes the proof of Lemma \ref{L:speed}.
\qed
\end{proof}

\medskip For $0\le x <1$ and $T>0$, let $\be_x^{T}$ denote the law
\begin{equation}
\be_x^{T}(\cdot)=\be(\cdot|Y_0=x; \sup_{0\le s \le T} Y_s \le 1)
\end{equation}
and let
\begin{equation}
\be^\infty_x=\lim_{T\to \infty} \be_x^T
\end{equation}
be the weak limit of $\be^T_x$, which may be described with Pinsky's result \cite{pinsky}.

\begin{lemma}
\label{L:int-1/2}
For any $\eps>0$, for any $\eta>0$, there exists $t_0$ such that for all $t\ge t_0$, and for all large enough $T>0$,
\begin{equation}
\be_{1/2}^T\left(\left|\frac1t\int_0^t Y_sds - m_0\right| >\eps\right) \le \eta.
\end{equation}
\end{lemma}

\begin{proof}
Let $t>0$. As $t\to \infty$, we know by the ergodic theorem for one-dimensional diffusions (Theorem V.53.1 in \cite{rw}), and the above calculations that, $\be_{1/2}^\infty-$almost surely,
\begin{equation}
\lim_{t\to \infty} \frac1t\int_0^t Y_s ds = m_0.
\end{equation}
Thus this convergence holds in $\be^\infty_{1/2}$-probability as well and we may choose $t_0$ large enough that
\begin{equation}
 \label{int1/2-1}
\be_{1/2}^\infty\left(\left|\frac1t \int_0^tY_s ds - m_0\right|>\eps\right) \le \eta/2
\end{equation}
for all $t\ge t_0$. Let us fix any $t\ge t_0$. Since $\be^T_{1/2}$ converges weakly towards $\be^\infty_{1/2}$, and since integration over the compact interval $[0,t]$ is a continuous functional, we conclude that
$$
\be^T_{1/2}\left(\int_0^t Y_s ds \in B \right) \to \be^\infty_{1/2}\left(\int_0^t Y_s ds \in B\right)
$$
for all Borel set $B\subset \mathbb{R}$, as $T\to \infty$. Taking $B=[(m_0-\eps)t,(m_0+\eps)t]$, we may choose $T_0$ large enough that for all $T\ge T_0$,
\begin{equation}\label{int1/2-2}
\left|\be_{1/2}^\infty\left(\left|\frac1t \int_0^tY_s ds - m_0\right|>\eps\right)- \be_{1/2}^T\left(\left|\frac1t \int_0^tY_s ds - m_0\right|>\eps\right)\right|\le \eta/2
\end{equation}
Combining (\ref{int1/2-1}) and (\ref{int1/2-2}) gives the result.
\qed
\end{proof}

\medskip The next step is to extend Lemma \ref{L:int-1/2} to a similar convergence type of result, but where the starting point $x$ is not necessarily equal to 1/2, while keeping the estimates uniform in $x$.

\begin{lemma}
\label{L:int-x}
For any $\eps>0$, for any $\eta>0$, there exists $t_0$ such that for all $t\ge t_0$, for all $x\in [0,1)$, and for all $T$ large enough,
\begin{equation}
\be_{x}^T\left(\left|\frac1t\int_0^t Y_sds - m_0\right| >\eps\right) \le \eta.
\end{equation}
\end{lemma}

\begin{proof}
We prove this by coupling. Consider two independent processes $Y^1$ and $Y^2$ sampled respectively from $\be^T_{1/2}$ and $\be^T_x$. Let $\tau=\tau(x,T)=\inf\{s>0: Y^1_s =Y^2_s\}$, and define
\begin{equation}
Y^3_s= Y^2_s\indic{s\le \tau} + Y^1_s \indic{s\ge \tau}
\end{equation}
It is easy to show that $Y^3$ has the same distribution as $Y^2$, i.e., its law is $\be^T_x$. Moreover, an application of Lemma \ref{Yf-ub} shows that the random variable $\tau$ is bounded above stochastically, uniformly in $T$ and $x\in[0,1)$. That is, for any $\eta$, there exists $t_1>0$ such that for all $T$ large enough and for all $x\in [0,1)$,
\begin{equation}\label{meet}
P(\tau>t_1)\le \eta.
\end{equation}
Indeed, the coupling time $\tau$ is smaller than the meeting time of two independent processes given by an unconditional square Bessel process of dimension 2 started at 1, with the diffusion $\be^\infty_0$. This meeting time is finite almost surely, which proves (\ref{meet}).
Let $\eps,\eta>0$. If we now choose $t$ large enough that $t_1/t \le \eps$ and $t>t_0$ from Lemma \ref{L:int-1/2}, we obtain:
\begin{align*}
\be_{x}^T\left(\left|\frac1t\int_0^t Y_sds - m_0\right| >2\eps\right) & \le  P(\tau > t_1) + P\left(\left|\frac1t\int_0^t Y^1_sds - m_0\right| >\eps\right) \\
& \le  \eta+ \be_{1/2}^x\left(\left|\frac1t\int_0^t Y_sds - m_0\right| >2\eps\right)
\end{align*}
Taking the limsup as $T\to \infty$, and using Lemma \ref{L:int-1/2}, we obtain
\begin{equation}\label{limsup}
\limsup_{T\to \infty} \be^T_x\left(\left|\frac1t\int_0^t Y_sds - m_0\right| >2\eps\right) \le 2 \eta
\end{equation}
for all $t\ge \max(t_0, t_1/\eps)$. Lemma \ref{L:int-x} is now easily deduced from (\ref{limsup}).
\qed
\end{proof}

\medskip Our next lemma shows that, given $\cE_a$, we are unlikely to spend a large amount of time below 0, and this amount can be controlled uniformly over $a$. In fact, the lemma states that once we reach a given level we are unlikely to spend more than a certain amount of time $z$ below it.

\begin{lemma} \label{L:tb}
 For any $\eps>0$, there exists $z>0$ such that for all $a>0$, and for all $0\le y <a$,
\begin{equation}
\bm_0\left( \left. \int_{\tau_y}^{\tau_a} \indic{X_s \le y} ds >z \right| \mathcal{E}_a \right) \le \eps
\end{equation}\label{max lev}
where $\tau_y=\inf\{s>0: X_s=y\}$. Similarly, there is $b>0$ such that for all $a>b$, for all $y \in (b,a)$,
\begin{equation}\label{minlev}
\bm_0\left( \inf_{\tau_y \le s \le \tau_a} X_t < y-b\right) \le \eta.
\end{equation}

\end{lemma}

\begin{proof}
For $s\ge 0$, let $\tilde X_s = X_{\tau_y+s}-y$ and let $\tilde L(s,w)= L(\tau_y+s,y+w)$. By the Markov property, it is easy to check that, given $\cE_a$, and given $\cF_{\tau_y}$, the process $\tilde X$ has the law $\bm_0(\cdot | \tilde \cE)$, where
\begin{equation}
\tilde \cE=\{ \tilde L(s,w) \le f(w), \text{ for all } 0 \le s \le \tilde \tau_{a-y}\}
\end{equation}
and where
\begin{equation}
f(w):= 1-\tilde L(0,w) = 1-L(\tau_y, y+w).
\end{equation}
For $s\ge 0$, let
\begin{equation}
Z_s = \tilde L(\tilde \tau_{a-y},-s)-\tilde L(0,-s)
\end{equation}
be the local time at level $y-s$ accumulated by $\tilde X$ after hitting $y$. Then note that by the occupation formula,
\begin{equation}
\int_{\tau_y}^{\tau_a} \indic{X_s \le y} ds = \int_0^{\tilde \tau_{y-a}} \indic{\tilde X_s \le 0}ds  =\int_0^\infty
 Z_s ds.
\end{equation}
By the Ray-Knight theorem, given $\tilde\cE$, and given $Z_0=x\in[0,1)$, $(Z_s,s\ge 0)$ has the law
\begin{equation}\label{Zcond}
\mathbb{Z}^f_x := \mathbb{Z}_x(\cdot |\{ Z_w \le  f(w), \text{ for all } w\ge 0\}),
\end{equation}
where $\mathbb{Z}_x$ denotes the law of a Bessel process of dimension 0 started from $Z_0=x$, i.e., $\mathbb{Z}_x$ is the Feller diffusion started from $x$.
 (Note that the event in the right-hand side of (\ref{Zcond}) is an event of positive probability for any given $x<f(0)$, since Feller diffusions become extinct almost surely). By Lemma \ref{Yf-ub} applied to the diffusion $Z$ rather than $Y$, for any $x<f(0)$, the conditional law $\mathbb{Z}^f_x$ is stochastically dominated by $\mathbb{Z}_x$. Using for instance the branching property of Feller diffusions, this is itself dominated by $\mathbb{Z}_1$, since $x < f(0) \le 1$. Thus, letting $\mu(dx)$ denote the law on $[0,f(0)]$ of $Z_0$,
\begin{align*}
\bm_0\left( \left. \int_{\tau_y}^{\tau_a} \indic{X_s \le y} ds >z \right| \cF_{\tau_y};\mathcal{E}_a \right) &= \bm_0\left( \left.
\int_0^{\tilde \tau_{y-a}} \indic{\tilde X_s \le 0}ds >z \right| \tilde \cE\right)\\
&= \int_0^1 \mu(dx) \bm_0\left( \left.\int_0^{\infty} Z_s ds  >z \right| Z_0=x; \tilde \cE \right)  \\
& \le  \int_0^1 \mu(dx) \mathbb{Z}^f_x\left( \int_0^\infty Z_s ds > z \right)\\
&\le  \int_0^1 \mu (dx) \mathbb{Z}_1 \left( \int_0^\infty Z_s ds > z \right)\\
&\le  \mathbb{Z}_1\left( \int_0^\infty Z_s ds > z \right).\\
\end{align*}
Now, under $\mathbb{Z}_1$, $(Z_s,s\ge 0)$ is almost surely continuous and becomes extinct in finite time, thus $\int_0^\infty Z_s ds<\infty$ almost surely, and the right-hand side in the above inequality can be made arbitrarily small for large enough $z$. Taking the expectation to average out the conditioning of $\cF_{\tau_y}$ finishes the proof of the first part of Lemma \ref{L:tb}. The second part (\ref{minlev}) also follows from the same method, details are left to the reader.
\qed\end{proof}

\medskip We now show how Lemma \ref{L:tb} can be applied to prove a first piece of the result in Theorem \ref{T:ballistic}: it is shown that if $y<a$ is given (we want to think of $y$ large but fixed, and $a \to \infty$), then given $\cE_a$ it has taken no more than approximately $m_0 y$ units of time to reach $y$.

\begin{lemma} \label{L:hit1}
 For any $\eps, \eta>0$ there exists $y_0$ large enough that if $y\ge y_0$, and for all sufficiently large $a>0$,
\begin{equation}
\bm_0(\tau_y > m_0y(1+\eps)|\mathcal{E}_a) \le \eta
\end{equation}
where $\tau_y =\inf\{s>0: X_s = y\}$.
\end{lemma}

\begin{proof}
We start by noticing that for any $z\ge 0$,
\begin{align*}
\bm_0(\tau_y > m_0y(1+\eps)|\mathcal{E}_a)  & \le  \bm_0\left( \left. \int_0^y L(\tau_y,w)dw > m_0 y(1+\eps)-z\right|\cE_a\right)\\
& \   + \bm_0\left( \left. \int_0^{\tau_a} \indic{X_s \le 0} ds > z\right| \cE_a\right). \\
\end{align*}
Thus if we choose $z$ as in Lemma \ref{L:tb} applied to $y=0$, we have, for any $a>0$, and for any $y\ge y_1:=2z/(m_0\eps)$
\begin{align}
\bm_0(\tau_y > m_0y(1+\eps)|\mathcal{E}_a) &\le \eta + \bm_0\left( \left. \int_0^y L(\tau_y,w)dw > m_0 y(1+\eps/2)\right|\cE_a\right)\nonumber \\
&\le  \eta + \bm_0\left( \left.\int_0^y L(\tau_a,w)dw > m_0 y(1+\eps/2)\right|\cE_a\right). \label{tb0}
\end{align}
For $w\ge 0$, let $Y_w= L(\tau_a, a-w)$. Under $\bm_0$, recall that by the Ray-Knight theorem, $(Y_w,w\ge0)$ is a strong Markov process which has the law of a square planar Bessel process for $(0\le w\le a)$ and a Feller diffusion for $w\ge a$. Now, conditionally on $\cE_a$, and conditionally on $Y_a= x \in(0,1)$, it follows easily from the strong Markov property at time $a$ that $(Y_s,0 \le s\le a)$ has the law of a square planar Bessel bridge conditioned on $\{\sup_{s\le a} Y_s \le 1\}$. That is, if we further condition on the position $Y_a=x$, the part of the constraint on $Y_w$ for $w\ge a$ becomes irrelevant.

\mn We now appeal to the following time-reversal argument: let $(Y_s,s\ge 0)$ be a square Bessel bridge of dimension 2 with $Y_0=0$ and $Y_a=x$, and let
\begin{equation}\label{rev0}
Y^\leftarrow_w=Y_{a-w}, \ \ 0\le w \le a
\end{equation}
be the time-reversed process. Then
$(Y^\leftarrow_w,0\le w \le a)$ is itself a square Bessel bridge of dimension 2 with $Y^\leftarrow_0=x$ and $Y^\leftarrow_a = 0$. (This follows quite easily from the rotational invariance of Brownian motion and from the fact that a Brownian bridge presents the same time-reversibility.) Furthermore, note that by Lemma \ref{Yf-ub}, a square Bessel bridge from $x$ to 0, conditioned on $\{Y_s \le 1 \text{ for all } s \le a \}$, can be related to the measure $\be_x^a$ in the following fashion:
\begin{align}
\be_x^a(\cdot | Y_a = 0) &= \lim_{\delta \to 0} \be_x^a (\cdot | Y_a \le \delta) \nonumber \\
&\preceq \be_x^a(\cdot )
\end{align}
where $\preceq$ stands for stochastic domination. Therefore, taking $Y_w = L(\tau_a,a-w)$,
\begin{align}
& \bm_0\left( \left. \int_0^y L(\tau_a,w)dw > m_0 y(1+\eps/2)\right|\cE_a, Y_a=x\right) \nonumber\\
&= \be_0^a(\int_{a-y}^a Y_w dw >m_0y(1+\eps/2) | Y_a = x)\nonumber \\
&= \be_x^a\left(\int_{0}^y Y_w dw >m_0y(1+\eps/2) | Y_a = 0\right)\nonumber \\
&\le  \be_x^a\left( \frac1y \int_0^y Y_w dw > m_0 (1+\eps/2)\right) \label{rev1}
\end{align}
By Lemma \ref{L:int-x}, we may choose $y_2$ large enough that if $y\ge y_2$ and for all large enough $a$, the right-hand side of (\ref{rev1}) is smaller than $\eta$. Thus for $y\ge y_1 \vee y_2$, and for all large enough $a$, we have by (\ref{tb0}) and unconditioning on the position $Y_a$ in (\ref{rev1}),
\begin{equation}
\bm_0(\tau_y > m_0y(1+\eps)|\mathcal{E}_a) \le 2\eta
\end{equation}
as required.
\qed\end{proof}

\medskip We now prove a bound in the other direction for the hitting times of certain levels. To start with, we need an a priori bound that says that it is unlikely for $L(\tau_a, 0)$ to be close to 1 when we condition on $\cE_a$.

\begin{lemma} \label{L:L(0)}
For any $\eta>0$, there is a $\delta >0$ such that
\begin{equation}\label{L(0)}
\bm_0(L(\tau_a, 0)\ge 1-\delta |\cE_a) \le \eta,
\end{equation}
for all large enough $a>0$.
\end{lemma}

\begin{proof}
By Lemma \ref{Yf-ub} and the Ray-Knight theorem, we observe that the random variable $L(\tau_a, 0)$, conditionally given $\cE_a$, is stochastically dominated by the squared modulus of a two-dimensional Brownian motion at time $a$, conditioned to be smaller than 1. However, the modulus at time $a$ is an exponential random variable with mean $\sqrt{a}$, so (\ref{L(0)}) follows easily.
\qed
\end{proof}

\begin{lemma} \label{L:hit2}
For any $\eps, \eta>0$ there exists $y_3$ large enough that if $y\ge y_3$, and for all sufficiently large $a>0$,
\begin{equation}
\bm_0(\tau_y < m_0y(1-\eps)|\mathcal{E}_a) \le \eta.
\end{equation}
\end{lemma}

\begin{proof}
The proof proceeds basically through the same steps as Lemma \ref{L:hit1}, but there are a few changes. Let $z$ be as in Lemma \ref{L:tb}, and let $2z/\eps=:y_4<y<a$. On the event $E(y,z)$ that $X$ doesn't spend more than $z$ units of time after $\tau_y$ below level $y$, we get:
\begin{align*}
E(y,z) \cap \{ \tau_y < m_0y(1-\eps)\} & = E(y,z) \cap \left\{\int_{-\infty}^y L(\tau_y,w)dw <m_0 y(1-\eps)\right\} \\
& \subset  \left\{\int_0^y L(\tau_a,w)dw <m_0 y(1-\eps)+z\right\}\\
& \subset  \left\{\int_0^y L(\tau_a,w)dw <m_0 y(1-\eps/2) \right\}.
\end{align*}
Define $(Y_w,w\ge 0)$ to be, as usual, $Y_w = L(\tau_a, a-w)$, for any $w\ge 0$. Recall that $Y$ is an inhomogeneous diffusion, or more precisely, a square Bessel process of dimension 2 on $[0,a]$, and a Feller diffusion on $[a, +\infty)$. Fix $\delta >0$ as in Lemma \ref{L:L(0)}, and note that by optional stopping, since $Z$ is a $\mathbb{Z}_{1-\delta}$-martingale,
$$
\mathbb{Z}_{1-\delta}\left(\sup_{s>0} Z_s <1\right) = \delta.
$$
Now, by Lemma \ref{L:int-x}, we can choose $y_0$ such that if $y\ge y_0$, for all $x \in (0,1)$ and all $b>0$ large enough,
\begin{equation}\label{eta delta}
\be_x^b\left( \frac1y \int_0^y Y_sds <m_0y(1-\eps/2)\right) \le {\eta}{\delta}
\end{equation}
Therefore,
\begin{align}
&\bm_0\left(\left. \int_0^y L(\tau_a,w)dw <m_0 y(1-\eps/2)\right|\cE_a\right)  \nonumber\\
\le &  \eta + \bm_0\left(\left.\int_0^y L(\tau_a,w)dw <m_0 y(1-\eps/2); Y_a \le 1-\delta\right|\cE_a\right) \nonumber\\
\le & \eta + \frac1\delta \be_0^a\left(\int_{a-y}^a Y_sds <m_0 y(1-\eps/2); Y_a \le 1-\delta\right)  \nonumber\\
\le & \eta + \frac1\delta \be_0^a\left(\int_{a-y}^a Y_sds <m_0 y(1-\eps/2)\right) \label{eta delta 2}
\end{align}
The idea is now to condition upon the position $Y_{a-y}=x$. Conditionally on this event,
$$
\be_0^a\left(\left.\int_{a-y}^a Y_sds <m_0 y(1-\eps/2)\right|Y_{a-y}=x\right) = \be_x^y\left( \frac1y \int_0^y Y_s ds \le m_0(1-\eps/2)\right).
$$
However, by Lemma \ref{Yf-ub}, $\be_0^y\preceq \be_0^b$ for any $b>y$. Thus
\begin{eqnarray*}
\be_0^a\left(\left.\int_{a-y}^a Y_sds <m_0 y(1-\eps/2)\right|Y_{a-y}=x\right) & \le & \be_x^b\left( \frac1y \int_0^y Y_s ds \le m_0(1-\eps/2)\right). \\
& \le & \eta \delta
\end{eqnarray*}
by our choice of $y\ge y_0$ and by taking $b$ sufficiently large that (\ref{eta delta}) holds. Plugging this into (\ref{eta delta 2}), we obtain:
$$
\bm_0\left(\left.\int_0^y L(\tau_a,w)dw <m_0 y(1-\eps/2)\right|\cE_a\right) \le \eta + \frac1{\delta} \eta \delta = 2 \eta.
$$
This completes the proof of Lemma \ref{L:hit2}.
\qed\end{proof}

\medskip We are now ready to finish the proof of Theorem \ref{T:ballistic}.

\mn \emph{Proof of Theorem \ref{T:ballistic}.} The proof of Theorem \ref{T:ballistic} is divided into two steps, a lower bound and an upper bound. We start with the lower bound. We want to show that for any $\eps,\eta>0$, there exists $t_3$ large enough that for all $t\ge t_3$, and for all $a>0$ sufficiently large,
\begin{equation}
\label{lb0}
\bm_0(X_t < \gamma_0t(1-\eps)|\cE_a) \le \eta.
\end{equation}
Indeed, if this holds, then it follows by weak convergence that any subsequential limit $\mathbb{Q}$ of $\bm_0(\cdot | \cE_a)$ satisfies: for any $\eps$, any $\eta>0$, there exists $t_3$ such that for all $t\ge t_3$
\begin{equation}\label{lb1}
\mathbb{Q}(X_t<\gamma_0t(1-\eps) )\le \eta
\end{equation}
because the canonical projection map $X \mapsto X_t$ is a continuous map for the Skorokhod topology. Hence
$\mathbb{Q}( {X_t}/t - \gamma_0 < -\eps) \le \eta$,
and we conclude:
$$
\mathbb{Q}( \frac{X_t}t - \gamma_0 < -\eps)  \xrightarrow[t\to \infty]{} 0
$$
which is, as claimed, the lower bound required for the proof of Theorem \ref{T:ballistic}.
Let us thus turn to (\ref{lb1}) and fix $\eps, \eta>0$ with $\eps<1$, and choose $y_0$ as in Lemma \ref{L:hit1}. For $t_4 = 4y_0m_0$, and for $t\ge t_4$, let $y= \gamma_0 t(1-\eps/2) \ge y_0$. Thus,  for all $a$ sufficiently large,
$$
\bm_0( \tau_y \le t(1-\eps/4)| \cE_a) \le \eta.
$$
Having reached level $y= \gamma_0t(1-\eps/2)$ by time $t(1-\eps/4)$, the only way $X_t$ can be below $\gamma_0(1-\eps)t$ is if $X$ reaches again $\gamma_0t(1-\eps)$ after time $\tau_y$. By (\ref{max lev}) in Lemma \ref{L:tb}, if $t\ge t_5 = 4b/\eps$ (where $z$ is as in Lemma \ref{L:tb}), then this occurs with probability at most $\eta$ for all large enough $a$. Thus we conclude, for $t\ge t_3:= t_4 \vee t_5$, for all large enough $a$,
\begin{equation}
\bm_0(X_t < \gamma_0t(1-\eps)|\cE_a) \le 2\eta.
\end{equation}
This concludes the proof of the lower bound. We now turn to the proof of the upper bound, where we wish to prove that for all $\eta, \eps>0$, there is $t_6$ large enough that for all $t\ge t_6$, and for all $a>0$ large enough,
\begin{equation}
\bm_0(X_t > \gamma_0t(1+\eps)|\cE_a) \le \eta.
\label{ub0}
\end{equation}
However, note that the event $\{X_t > \gamma_0t(1+\eps)\}$ is contained in the event $\{\tau_{y} \le t\}$ where $y = \gamma_0t (1+\eps)$. By Lemma \ref{L:hit2}, if $y\ge y_3$, in particular if $t\ge t_6:=y_3 m_0$, then it follows:
\begin{align*}
 \bm_0(X_t > \gamma_0t(1+\eps)|\cE_a)  & \le  \bm_0(\tau_y \le t |\cE_a) \\
& \le  \bm(\tau_y \le \frac{ym_0}{1+\eps} | \cE_a)\\
& \le \eta,
\end{align*}
as desired. This completes the proof of Theorem \ref{T:ballistic}. \qed

\section{Random walk with bounded local time}

Throughout this section we assume
\begin{equation}\label{new5abc}
L_0 \ge 2.
\end{equation}

\medskip We need to introduce some notation. Let
\begin{equation} \label{tau_k}
\tau_k := \inf\{i:S_i = k\}
\end{equation}
be the first hitting time of $k \ge 0$. We then define $\cB_k,
\cB_k^+$ to be the events
\begin{equation} \label{B_k}
\cB_k := \cA_{\tau_k} = \{L(\tau_k,x, \om) \le L_0 \text{ for all
}x\}
\end{equation}
and
\begin{equation}\label{B+_k}
\cB_n^+ := \cB_n \cap \{S_i > 0 \text{ for } 1\le  i \le \tau_n\}.
\end{equation}
Thus $\cB_n^+$ occurs if the sample path minus its endpoints stays
strictly between its initial point at 0 and its final point at
$n$. Thus the maximum value of the points is $n$ and this is taken
on for the first time at the endpoint and
necessarily, the length of the path equals $\tau_n$.
Moreover, the sample path through time $\tau_k$  visits each value
$x$ at most $L_0$ times. The event $\cB^+_k$ will play a major role in our analysis, since it can be interpreted as having a regenerating level immediately at the starting point. We shall make use of the following $\si$-fields:
$$
\cF_n = \si \{S_i, i \le n\}, \quad \cF_\infty = \bigvee_{n \ge 0}
\cF_n,
$$
and
$$
\cG_k = \cF_{\tau_k}.
$$

\begin{lemma} \label{L:Bk-vs-Bk+} There exists some constant $C_3 > 0$ such that
\begin{equation}\label{6.1}
P(\cB_k^+) \ge C_3 P(\cB_k), k \ge 1.
\end{equation}
\end{lemma}
\begin{proof} Fix $k$ and let $\rho$ be the last time before
$\tau_k$ at which the random walk visits 0, i.e.,
$$
\rho = \max\{i < \tau_k: S_i = 0\}.
$$
Note that $S_{\tau_k} = k > 0$ for $k \ge 1$. Therefore, $S_i > 0$
for $\rho < i \le \tau_k$. Consequently, a decomposition with
respect to the value of $\rho$ shows that
\begin{align}
P(\cB_k) &= \sum_{j=0}^{\infty} P(\rho = j, \cB_k)\nonumber\\
&\le\sum_{j=0}^{\infty} P(S_j = 0, L(j,x) \le L_0 \text{ for
all
}x,\tau_k > j \nonumber\\
& \ \ \  \text{ and } S_n - S_j > 0 \text{ for } 1 \le n-j \le \tau_k
-j, L(\tau_k,x) -
L(j,x) \le L_0 \text{ for all } x)\nonumber\\
&=\sum_{j=0}^{\infty} P(S_j = 0, L(j,x) \le L_0 \text{ for all
}x, j < \tau_k)P(\cB_k^+) \label{new8}
\end{align}
But for any $x$, on the event $\{S_n=y, L(n,x) \le L_0 \text{ for
all }x\}$ it holds
\begin{align*}
&P(L(n+2L_0+2,y) \ge L_0+1\big|S_0, \dots, S_n) \\
& \ \ \ge P(S_{n+2i+1} = y+1, S_{n+2i+2} = S_n = y \text{ for } 0 \le i \le L_0)\\
&\ \ \ge 2^{-L_0-1} > 0.
\end{align*}
It follows easily from this that
\begin{equation}\label{new7}
P(L(j,x) \le L_0 \text{ for all }x) = P(\cA_j) \le
C_4e^{-C_5j}
\end{equation}
for some constants $0 < C_i < \infty$. In turn, this implies
$$
\sum_{j=0}^{\infty} P(S_j = 0, L(j,x) \le L_0 \text{ for all
}x, j < \tau_k)\le \sum_{j=0}^\infty P(L(j,x) \le L_0 \text{ for all }x) < \infty,
$$
so that (\ref{6.1}) follows from (\ref{new8}). \qed
\end{proof}

\medskip We need sharper information about possible weak limits of
$P(\cdot|\cB_r)$. This will be given in the following lemma. We
define
\begin{equation}\label{new7.1}
\cC_k:= \{S_i > k \text{ for all }i > \tau_k\}.
\end{equation}

\noindent {\bf Remark 1.} We are going to study weak limit points
of the measures $P(\cdot \big|\cB_r)$ as $r \to \infty$. Note
that each $\tau_n < \infty$ a.s. $[P]$, so conditioning on $\cB_r$
is the same as conditioning on $\cB_r \cap \{\tau_n < \infty\}$
for any $n$, including $n=r$, possibly. This does not
automatically say that for a limit point $Q$ of $P(\cdot
\big|\cB_{r_i})$ it holds $Q(\tau_n < \infty)= 1$ for all $n$.
In fact this will be false for $n < 0$. But it is correct for $n
\ge 0$. Indeed, the case $n=0$ is trivial, since $\tau_0 =0$ a.s.
$[P]$. For $r > n> 0$,
\begin{equation}\label{hit0}
P(\tau_n > t, \cB_r) = E\left( \indic{\tau_n > t}
P(\cB_r\big|\cF_t)\right) \le P(\tau_n > t)P(\cB_{r-n}).
\end{equation}
To see this, note that if the walk is at a position $m < n$ at time
$t$, then for $\cB_r$ to occur the local time has to be $\le L_0$
as the walk moves from $m$ to $r$, which is an interval of length
at least $r-m \ge r-n$: this implies (\ref{hit0}). Therefore, by (\ref{new6.5z}) below,
$$
P(\tau_n > t\big| \cB_r)\le P(\tau_n > t)\frac
{P(\cB_{r-n})}{P(\cB_r)}  \le P(\tau_n > t) 2^n.
$$
For fixed $n > 0$ we can make the limsup of the right hand side
here as $r \to \infty$ as small as we like by taking $t$ large.
Thus $Q(\tau_n = \infty) = 0$ for each $n > 0$.

The following lemma is the first of two crucial steps in the proof of Theorem \ref{T:rw}.

\begin{lemma} \label{L:crucial1} There exists a constant $0 \le C_4 < \infty$ such that
\begin{equation}\label{new5z}
\limt [P(\cB_t)]^{1/t} = e^{-C_4}
\end{equation}
and for all $t \ge 0$,
\begin{equation}\label{new6z}
P(\cB_t) \ge e^{-C_4t}.
\end{equation}
In addition, for all $s,t \ge 0$:
\begin{equation}\label{new6.5z}
P(\cB_t) \le 2^s P(\cB_{t+s})
\end{equation}
Further,
\begin{equation}\label{new7z}
 P(\cB_n^+) \sim
C_6 e^{-C_4n}
\end{equation}
for a suitable constant $C_6 > 0$.
\end{lemma}

\begin{proof} For (\ref{new5z}) and (\ref{new6z}),
we merely have to observe
that
\begin{equation}\label{new7.75z}
P(\cB_{s+t}) = P(\cA_{\tau_{s+t}}) \le
P(\cA_{\tau_s})P(\cA_{\tau_t}) = P(\cB_s)P(\cB_t),
\end{equation}
because if the random walk $\{S_n\}$ reaches the level $s+t$ at
time $\tau_{s+t}$, with $\sup_x L(\tau_{s+t},x) \le L_0$, then the
random walk must first reach $s$ at time $\tau_s$ with $\sup_x
L(\tau_s,x) \le L_0$ and then the random walk starting at $s$ must
reach $s+t$ with $\sup_x [L(\tau_{s+t},x)-L(\tau_s,x)] \le L_0$. Thus $P(\cB_t)$ forms a submultiplicative sequence, and it follows that $\lim_{t\to \infty} P(\cB_t)^{1/t} = e^{-C_4}$ exists. It is obvious that $C_4 \ge 0$, and from (\ref{new6.5z}), proved below, we get that $C_4\le \log(2) < \infty$. Moreover it is well-known that by submultiplicativity, $-C_4 = \inf_{t\ge 1}\{ \log P(\cB_t)/t\}$, hence $P(\cB_t) \ge e^{-C_4 t}$ for all $t\ge 1$.

As for (\ref{new6.5z}), this follows from the simple fact that (by
definition) the random walk arrives at $t$ for the first time at
$\tau_t$, so that $S_{\tau_t} =t$. If then the random walk takes
one step to the right it arrives for the first time at $t+1$ at
time $\tau_t +1$. Moreover, $\sup_x L(\tau_t +1,x) \le 1 \lor
\sup_x L(\tau_t,x)$, because the random walk visits a new point at
$\tau_t+1$. Thus, if $\cB_t$ occurred, then also $\cB_{t+1}$ occurs
in this case. Hence
$$
P(\cB_{t+1}) \ge P(\cB_t)P(S_{\tau_t+1} = S_{\tau_t} + 1) =
\frac 1 2 P(\cB_t).
$$
Induction on $s$ now yields (\ref{new6.5z}).

The proof of (\ref{new7z}) is much more involved. However, it is
closely related to Lemma 2 in Kesten \cite{kesten}. In analogy with the $L_n$ from this
reference we introduce the further event $\cL_n$ which is roughly
speaking the event that $\cB^+_n$ occurs (so that 0 is a regeneration level) and there is no other regeneration level between 0 and $n$. To give the formal definition, we define the shift $T_n$ by
$$
(T_n\om)_j = \om_{\tau_n+j}
$$
We then take $\cB^+_0$ to be the certain event, $\cL_0$ the empty
event,  and $\cL_1 = \cB^+_1$ the event $\{S_0 =0, S_1 = 1\}$.
Further, for $n \ge 2$
\begin{equation}\label{new10z}
\cL_n := \cB^+_n \cap \{ \forall k <n, \ T_k \om \notin \cB^+_{n-k}\}.
\end{equation}
The last property says that a sample path $(\om_0, \om_1, \dots,
\om_m)$ in $\cL_n$ cannot be decomposed into two pieces $(\om_0,
\dots, \om_j)$ and $(\om_j, \dots, \om_m)$ with the first part
minus its endpoint lying strictly to the left of  $\om_j$ and the
second part lying strictly to the right of $\om_j$ (except for its
initial point). The first part in such a decomposition would
belong to $\cB_j^+$ and the second part would be a translate of a
path in $\cB_{n-j}^+$.



Of course $\{T_k \om \in \cB^+_{n-k}\}$ is the event that
$\cB_{n-k}^+$ occurs for the shifted sequence $T_k\om =
(\om_{\tau_k}, \om_{\tau_k+1}, \dots)$.  Since $\cB_{n-k}^+$
depends only on $(\om_0, \dots,\om_{\tau_{n-k}})$ we shall
occasionally abuse notation and write $(\om_{\tau_k},
\om_{\tau_k+1}, \dots, \om_{\tau_n}) \in \cB^+_{n-k}$ instead of
$T_k \om \in \cB_{n-k}^+$.

The main step will be to show that
\begin{equation}\label{new11z}
P(\cB^+_n) = \sum_{j=1}^n P(\cL_j)P(\cB_{n-j}^+), \quad n\ge
1.
\end{equation}
This relation holds by convention if $n=1$, so assume $n \ge 2$
and that $\cB^+_n$ occurs. Then define $k$ to be minimal so that
$\cB_k^+ \cap\{T_k\om \in \cB_{n-k}^+\}$ occurs. This minimal
index is well defined because the event $\cB_n^+ \cap\{T_n\om \in
\cB_0^+\} = \cB_n^+$ occurs. Of course the minimal index is
unique. We claim that for this minimal $k$ the event $\cL_k$
occurs. Indeed, note that $\cB_k^+$ occurs, so that by the definition (\ref{new10z}) with $n$ and $k$ replaced by $k$ and $j$, if $\cL_k$ fails, then it must be that $\{\forall j<k, T_j \omega \notin \cB^+_{k-j}\}$ fails, i.e., there is $j<k$ such that $T_j \omega \in \cB^+_{k-j}$. Since $\omega \in \cB^+_n$, this implies that $T_j \omega \in \cB^+_{n-j}$ as well, and it is obvious that $\cB^+_j$ must hold as well since $\cB^+_n$ holds. This contradicts the minimality of $k$, and hence $\cL_k$ holds. Since $T_k \omega \in \cB^+_{n-k}$ by definition, it follows immediately that
$$
P(\cB_n^+) \le \sum_{k=1}^n P\big(\cL_k \cap\{T_k\om \in
\cB_{n-k}^+\}\big).
$$
But $\cL_k \in \cG_k$, because the occurrence of $\cL_k$ depends
on $(\om_0, \dots, \om_{\tau_k})$ only. (Recall that $\cG_k = \cF_{\tau_k}$ by definition). Thus, by the strong Markov
property
\begin{equation}\label{new14z}
P(\cB_n^+) \le \sum_{k=1}^n P(\cL_k)P(T_k\om \in
\cB_{n-k}^+) = \sum_{k=1}^n P(\cL_k) P(\cB_{n-k}^+).
\end{equation}
To prove the opposite inequality fix a $k \in \{1,\ldots, n\}$ and assume
the following two events occur:
\begin{equation}\label{new15z}
\cL_k \text{ and } \om':= T_k\om = (\om_{\tau_k}, \om_{\tau_k+1},
\dots) \in \cB_{n-k}^+.
\end{equation}
Then $\om$ is such that
\begin{equation}\label{new16z}
1 \le \om_\ell \le k-1 \text{ for }1 \le \ell < \tau_k,
\om_{\tau_k} = k,
\end{equation}
and
\begin{equation}\label{new17z}
k+1 \le \om_{\tau_k+ \ell} = \om'_\ell \le n-1 \text{ for } 1 \le
\ell \le \tau_n -1.
\end{equation}
Moreover, if $\tau'_{n-k}$ denotes the first hitting time of $n-k$
by the path $\om'$, then
\begin{align}
\sup_x L(\tau_n,x,\om) &\le \sup_x
L(\tau_k,x, \om) \lor \sup_x
L(\tau'_{n-k},x, \om') \lor 1 \nonumber\\
&=\sup_x L(\tau_k,x, \om) \lor \sup_x [L(\tau, n,x,
\om)-L(k,x,\om)] \lor 1 \nonumber \\
& \le L_0. \label{new18zy}
\end{align}
Together these
properties show that $\om \in \cB_n^+$. Thus the sample sequences
for which the events in (\ref{new15z}) occur
contribute $P(\cL_k)P(\cB_{n-k}^+)$ to $P(\cB_n^+)$. In
order to prove
\begin{equation}\label{new19z}
P(\cB_n^+) \ge \sum_{j=1}^n P(\cL_j)P(\cB_{n-j}^+)
\end{equation}
we therefore merely have to show that  (\ref{new15z}) can occur
only for one $k$. To see that this is indeed the case assume that
in addition to (\ref{new15z}) also
\begin{equation}\label{new20z}
\cL_j \text{ and } \om'':= T_j\om = (\om_{\tau_j}, \om_{\tau_j+1},
\dots) \in \cB_{n-j}^+
\end{equation}
occurs for some $j \ne k, j \in [1,n]$. For the sake of argument
let $j <k$. But then, we have on the one hand that $\cB_j^+$ occurs (by definition of $\cL_j$ or since $\cB_k^+$ occurs) and on the other hand, we also have that
\begin{equation}\label{new21z}
T_j\om  \in \cB_{k-j}^+ \text{ occurs}.
\end{equation}
But this contradicts the definition of $\cL_k$, so that (\ref{new15z}) and (\ref{new20z})
cannot hold simultaneously. This, in turn, implies (\ref{new19z})
and then finally (\ref{new11z}).

We can finally start on the proof of (\ref{new7z}) proper. Define
$$
f_n = e^{C_4n}P(\cL_n)  \text{ and } u_n = e^{C_4n}P(\cB_n^+).
$$
By our conventions just before (\ref{new10z})
$$
u_0 = 1, f_0 =0, u_1 = f_1 =(1/2) e^{C_4}.
$$
Moreover, by (\ref{new11z}) these quantities satisfy the renewal
equation
$$
u_n = \sum_{j=1}^n f_ju_{n-j}, \quad n \ge 1.
$$
In addition, by Lemma \ref{L:Bk-vs-Bk+} and (\ref{new6z})
$$
u_n = e^{C_4n} P(\cB_n^+) \ge  C_3e^{C_4n}P(\cB_n) \ge C_3 > 0
$$
is bounded away from 0, and $\limn [u_n]^{1/n} = 1$. By the renewal
theorem (see, e.g., Feller \cite[Theorems 2 and 3 in  12.3]{feller}), these facts imply
$$
\sum_{j=1}^\infty f_j = 1 \text{ and } \limn u_n = \frac 1\mu,
$$
where
\begin{equation}\label{new14z}
0 < \mu = \sum_{j=1}^\infty nf_n < \infty.
\end{equation}
Thus,
\begin{equation}\label{new14.5z}
P(\cB_n^+) \sim \frac 1{\mu} e^{-C_4 n}
\end{equation}
which proves (\ref{new7z}). The finishes the proof of Lemma \ref{L:crucial1}. \qed

\medskip We now move on to the second crucial step in the proof of Theorem \ref{T:rw}.

\begin{lemma} \label{L:crucial2}
We have:
\begin{equation}\label{new6.6z}
C_5 := \limn \frac {P(\cB_n^+)}{P(\cB_n)} \text{ exists and
}C_5 \ge C_3 > 0.
\end{equation}
Also for $\cE$ an event in $\cG_k$,
\begin{equation}\label{new6.7z}
\limn P(\cE \cap \cD_{k,n}\big|\cB_n) = C_5 e^{C_4k}P\left(\cE,\sup_x
L(\tau_k,x) \le L_0\right),
\end{equation}
where $\cD_{k,n}:=\{S_i > k \text{ for all }  \tau_k<i \le \tau_n\}$.
Finally,
\begin{equation}\label{P(Bn)}
\limn e^{C_4n}P(\cB_n) = C_7
\end{equation}
exists, where $0<C_7<\infty$.
\end{lemma}

\begin{proof} Let $\cR = \{0 < r_1 < r_2 \dots\}$ be a subsequence
along which the weak limit of $P(\cdot |\cA_{\tau_k})$ exists,
and let $Q(\cdot)$ be the value of this limit. The limit along
the subsequence $\cR$ will be denoted as $\lim_{r \in \cR}$.
Without loss of generality we may assume that also $\lim_{r \in
\cR} P(\cB_r^+)/P(\cB_r)$ exists (since it is a bounded sequence) and is at least $C_3$ (by
Lemma \ref{L:Bk-vs-Bk+}). (Later on we will prove that this limit does not depend on $\cR$ and hence $\lim_{r\to \infty} P( \cB^+_r) /P(\cB_r)$ exists.) Now let $\cE \in \cG_k$. Then
\begin{equation}\label{new15.1z}
Q(\cE \cap \cC_k)= Q(\cE, S_i > k \text{ for all }i > \tau_k)
= \lim_{N \to \infty} Q(\cE, S_i > k,  \tau_k < i \le
\tau_{k+N}).
\end{equation}
We want to show that this equals
\begin{equation}\label{new15zz}
e^{C_4k}P(\cE,\sup_x L(\tau_k,x) \le L_0) \lim_{r \in \cR} \frac
{P(\cB_r^+)}{P(\cB_r)}.
\end{equation}
To this end observe first that
$$
Q(\cE, S_i > k \text{  for } \tau_k < i \le \tau_{k+N})
= \lim_{r \in \cR} \displaystyle \frac{P(\cE, S_i > k \text{ for } \tau_k < i \le
\tau_{k+N}, \cB_r)}{P(\cB_r)},
$$
and secondly that for $r \ge k+N$ (because $S_{\tau_k} = k$)
\begin{align*}
&  \big|P(\cE, S_i > k
 \text{ for } \tau_k < i \le \tau_{k+N},
\cB_r)  - P(\cE, S_i > k \text{ for } \tau_k < i \le \tau_r,
\cB_r) \big|\\
&\le P(S_i = k \text{ for some }\tau_{k +N} < i \le \tau_r,
\sup_x L(\tau_r,x)\le L_0)\\
&\le P(\text{there exists some $\tau_{k+N} < i \le \tau_r$ for
which
$S_i = k$} \\
& \ \ \ \ \ \ \text{and $\sup_x L(i,x) \le L_0$ as well as }
\sup_x [L(\tau_r,x)-L(i,x)] \le L_0) \\
&\le P(S_i = k \text{ and $\sup_x L(i,x) \le L_0$ for some
$ \tau_{k+N}< i \le \tau_r$}) P(\cB_{r-k})\\
&\le P(S_i = k \text{ and $\sup_x L(i,x) \le L_0$ for some
$ \tau_{k+N}< i \le \tau_r$}) 2^k P(\cB_r)\\ 
&\le 2^kP(\sup_x L(\tau_{k+N},x) \le L_0)P(\cB_r)\\
&= 2^kP(\cB_{k+N})P(\cB_r) \le C_7 2^k e^{-C_4(k+N)}
P(\cB_r),
\end{align*}
for some constant $C_7$ independent of $k,r$ (use Lemma \ref{L:Bk-vs-Bk+}
and (\ref{new14.5z}) for the last inequality). Consequently, using (\ref{new15.1z}),
\begin{equation}
Q(\cE, S_i > k \text{ for all }i > \tau_k)
= \lim_{r \in \cR} \displaystyle \frac{P(\cE, S_i > k,
\tau_k < i \le \tau_r, \cB_r )}{P(\cB_r)}.
\label{new16zz}
\end{equation}
But if $S_i > k$ for $\tau_k < i \le \tau_r$, then
$$
L(\tau_r,x) = \begin{cases} L(\tau_k,x) &\text{ if } x \le k\\
L(\tau_r, x)- L(\tau_k,x) &\text{ if } x > k.
\end{cases}
$$
Therefore (use $\cE \in \cG_{k}$ and again $S_{\tau_k} = k$)
\begin{align}
P(\cE, S_i > k, \tau_k < i \le \tau_r,
\cB_r)
&= P(\cE, \cB_k)P(\cB_{r-k}^+).
\label{new16.5z}
\end{align}
Together with (\ref{new16zz}) and (\ref{new14.5z}) this proves the
desired (\ref{new15zz}).

We next claim that there exist events $\cE_k \in \cG_k$ such that
\begin{align}
\{\cE_k, S_i > k \text{ for all } i > \tau_k\}
& = \cE_k \cap \cC_k\nonumber\\
& =\{k \text{ is the smallest value of $n$ for which $\cC_n$
occurs}\}.\label{new17zy}
\end{align}
To see this, recall the definition for of $\cD_{j,k}$ for $j < k$,
$$
\cD_{j,k} = \{S_i > j \text{ for }\tau_j < i \le \tau_k\}.
$$
Then $\cC_j \cap \cC_k = \cD_{j,k} \cap \cC_k$ and consequently
$$
\cup_{0 \le j < k}( \cC_j \cap \cC_k) = \cC_k \cap \big[\cup_{0
\le j < k}\cD_{j,k}\big].
$$
The right hand side of (\ref{new17zy}) equals
$$
\cC_k \setminus \cup_{0\le j < k} (\cC_j \cap \cC_k) = \cC_k \cap
[\cup_{0 \le j <k}\cD_{j,k}]^c.
$$
This gives us (\ref{new17zy}) with $\cE_k$ equal to the complement
of $\cup_{0 \le j <k}\cD_{j,k}$.

\medskip We can now apply (\ref{new15zz}) with $\cE$ taken equal to $\cE_k$,
with the result that
\begin{align}
&   Q(k\text{ is the smallest value of $n$ for which $\cC_n$
occurs})
 =  Q(\cE_k \cap \cC_k)\nonumber \\
&    =e^{C_4k}P(\cE_k,\sup_x L(\tau_k,x) \le L_0)
\lim_{r \in \cR} \frac{P(\cB_r^+)}{P(\cB_r)}.
\label{new18zz}
\end{align}

Finally we shall show that
\begin{align}
&\sum_{k=0}^\infty Q(k\text{ is the smallest value of $n$
for which $\cC_n$ occurs})\nonumber \\
&=Q(\cC_k \text{ occurs for some $k \ge 0$}) =1.
\label{new15zy}
\end{align}
From this and (\ref{new18zz}) we can conclude that
\begin{equation}\label{new16zy}
C_5 :=\lim_{r \in \cR}\frac {P(\cB_r^+)}{P(\cB_r)}
\end{equation}
exists, is independent of $\cR$, and $\ge C_3$ by virtue of
Lemma \ref{L:Bk-vs-Bk+}. In view of (\ref{new16zz}) this will also show that
for all $\cE \in \cG_k$ the full limit
$$
\limn P(\cE \cap \cD_{k,n}\big | \cB_n) = C_5
e^{C_4k}P(\cE_k,\sup_x L(\tau_k,x) \le L_0)
$$
exists, and has the value given in (\ref{new6.7z}). Also (\ref{P(Bn)}) follows from (\ref{new7z}) and
(\ref{new16zy}).

It remains to prove (\ref{new15zy}). To this end we want to show
that $Q(\cC_k|\cG_k)$ is bounded from below. To prove this we
note that for each fixed $k$, any element of $\cG_k$ is up to
$Q$-null sets a finite or countable disjoint union of sets of the
form
$$
\cH(\boldsymbol \eta) = \{S_i = \om_i = \eta_i, 0 \le i \le m\},
$$
where $m < \infty$ and $\boldsymbol \eta = (\eta_0, \dots,
\eta_m)$ runs over the sequences which satisfy
\begin{equation}\label{new21zy}
\eta_0 = 0, \eta_i -\eta_{i-1} = \pm 1, 1 \le i \le m, \eta_j <
\eta_m = k, 0 < j < m.
\end{equation}
(Note that the requirements on $\eta$ in (\ref{new21zy}) are such
that $\tau_k = m$ for any sample point with $(S_0, \dots S_m) =
\boldsymbol \eta$. We can restrict ourselves to finite $m$,
because $Q(\tau_k = \infty) = 0$ by Remark 1.) Now, as before,
for any such $\boldsymbol \eta$
\begin{equation}\label{new20zy}
\frac{Q(\cH(\boldsymbol \eta), \cC_k)} {Q(\cH(\boldsymbol
\eta))}
 = \lim_{N \to \infty} \lim_{r \in \cR}
\frac {P(\cH(\boldsymbol \eta), \cD_{k,k+N}, \sup_x L(\tau_r,x) \le L_0)}{P(\cH(\boldsymbol
\eta), \sup_x L(\tau_r,x) \le L_0)}.
\end{equation}
This time we use that the denominator in the right hand side here
is at most
$$
P(\cH(\boldsymbol \eta), \sup_x L(\tau_k,x) \le
L_0)P(\cB_{r-k})
$$
(compare (\ref{new16.5z})). As in the lines following
(\ref{new15zz}) the numerator in the right hand side of
(\ref{new20zy}) is bounded below by
\begin{align*}
& P(\cH(\boldsymbol \eta), \sup_x L(\tau_k,x) \le L_0)
P(\cD_{k,k+N},
\sup_x [L(\tau_r,x) - L(\tau_k,x)] \le L_0)\\
& = P(\cH(\boldsymbol \eta), \sup_x L(\tau_k,x) \le L_0)
P(\cB_{r-k}^+)\\
&\ge C_3P(\cH(\boldsymbol \eta), \sup_x L(\tau_k,x) \le L_0)
P(\cB_{r-k}).
\end{align*}
It follows from these estimates that
$$
\frac{Q(\cH(\boldsymbol \eta), \cC_k)} {Q(\cH(\boldsymbol
\eta))} \ge C_3.
$$
Since this holds for all atoms $\eta$ of $\cG_k$ we conclude that
\begin{equation}\label{new21.2z}
Q(\cC_k|\cG_k) \ge C_3.
\end{equation}
The relation (\ref{new15zy}) is a simple consequence of
(\ref{new21.2z}) and the martingale convergence theorem. Indeed,
set
$$
Y_N = \indic{\cC_k \text{ occurs for
some $k \ge N$}}.
$$
Then, if we write $E^Q$ for expectation with
respect to $Q$, we have for each fixed $N$
$$
\lim_{k \to \infty} E^Q(Y_N \big |\cG_k) = Y_N \text{ a.s. }[Q].
$$
On  the other hand, for $k \ge N$, $E^Q(Y_N\big|\cG_k) \ge
Q(\cC_k\big|\cG_k) \ge C_3$, from which we deduce that
$$
Y_N \ge C_3 \ \ a.s.\  [Q].
$$  Thus $Q(Y_N =1)=1$ and
(\ref{new15zy}) holds. This finishes the proof of Lemma \ref{L:crucial2}. \qed
\end{proof}

\begin{lemma}\label{L:conv}
The sequence of measures $P(\cdot | \cB_n)$ converges weakly to a limit measure $Q$.
\end{lemma}

\begin{proof} Because the walk is nearest-neighbour, it is always the case that $P(\cdot | \cB_n)$ is tight: it thus suffices to prove uniqueness of the weak subsequential limits. Thus, let $Q$ be a weak limit along the subsequence $\cR$. Let $\nu_0=0$ and let $0 \le \nu_1 <
\nu_2 < \dots$ be the successive values of $\nu$ for which
$\cC_\nu$ occurs. (\ref{new15zy}) shows that $\nu_1 < \infty$ a.s.
$[Q]$, but the proof of (\ref{new15zy}) shows that all $\nu_i$ are
a.s. $[Q]$ finite. From this we will see that $Q(\cE) = \limn
P(\cE\big|\cB_n)$, with the limit taken along the sequence of
all integers, for any cylinder set $\cE$. Indeed, let $\cE \in
\cF_t$. Since $\tau_t \ge t$ (because $|S_{i+1} - S_i| \le 1$), we
have $\cF_t \subset \cG_t$, and so $\cE \in \cG_t$. Now let $\rho$
be the first $\nu_i \ge t$. Then $\cE = \cup_{s\ge t}[\cE
\cap\{\rho = s\}]$ and $\cE \cap \{\rho = s\} = \cE_s \cap \cC_s$
for some $\cE_s \in \cG_s$ (as in (\ref{new17zy})). Consequently,
$$
\limn P(\cE \cap\{\rho = s\}\big|\cB_n) \text{ exists} \text {
(by (\ref{new6.7z}))}.
$$
Also,
$$
\big|P(\cE\big|\cB_n) -
 \sum_{s=t}^{s=t+N} P(\cE \cap \{\rho= s\}\big|\cB_n)\big|
\le P(\rho > t+N\big|\cB_n).
$$
Finally,
\begin{equation}\label{new22z}
\lim_{N \to \infty} \limsup_{n \to \infty}P( \rho >
t+N\big|\cB_n) = 0,
\end{equation}
because if this fails, then (by the monotonicity in $N$) there
exists a sequence $\cR = \{r_1 < r_2 < \dots\}$ and an $\ep > 0$
such that
$$
Q(\rho \le t+N) \le 1-\ep \text{ for all } N,
$$
where $Q$ is the weak limit of $P(\cdot\big|\cB_{r_i})$. But we
have just seen that $\rho < \infty$ a.s. $[Q]$, so that
(\ref{new22z}) must hold. But then
$$
\limn P(\cE\big|\cB_n) = \sum_{s=t}^\infty \limn P(\cE \cap
\{\rho = s\}\big|\cB_n).
$$
This proves lemma \ref{L:conv}.\qed
\end{proof}

\bigskip
From now on $Q$ will be the (weak) limit of the probability
measures $P(\cdot \big|\cB_n)$ on $\Om$. Since $S_0 = 0$ and
$S_{n+1} - S_n = \pm 1$ with $P$-probability 1, it is also the
case that
\begin{equation}\label{new9}
Q(S_0 = 0) =1\text{ and } Q(S_{n+1} - S_n = \pm 1) =1.
\end{equation}
Also
\begin{equation}\label{new11}
Q(S_n = y \text{ for more than $L_0$ values of }n) = 0,
\end{equation}
because  for each fixed $n$ and all $r > n$,
$$
P((L_0+1)\text{-th
visit of $S$ to $y$ is at time $n$}\big|\cB_r) = 0.
$$
We remind
the reader that $\cC_k$ is defined in (\ref{new7.1}). We now come
to our main result, which describes the structure of $Q$ and is a more precise statement than Theorem \ref{T:rw}. Define
$\si_0=0$,
$$
\si_1 :=\inf\{\tau_\ell: \cC_\ell \text{ occurs}\},
$$
$$
\si_{j+1}:= \inf\{\tau_\ell > \si_j :\cC_\ell \text{ occurs}\},
$$
and, in agreement with Remark 2, let $\nu_j$ to be the  unique
value of $\nu$ for which $\si_j = \tau_\nu$. That is, $\sigma_j$ is the time at which the $j^{\text{th}}$ regeneration level occurs. Thus, by definition,
$$
S_n < \nu_j \text{ for } n < \tau_{\nu_j}, S_{\tau_{\nu_j}} =
\nu_j
$$
but
$$
S_n > S_{\si_j} =S_{\tau_{\nu_j}} = \nu_j \text{ for } n >
\tau_{\nu_j}.
$$
Moreover, if $n = \tau_s$ but $s$ is not one of the $\nu_j$, then
$S_t \le s$ for some $t > \tau_s$. Roughly speaking, the $\tau_j$
are the strict upward ladder epochs for the random walk $\{S_n\}$.
The $\si_j$ are special ladder epochs which make them into
regeneration times (in a sense to be made precise in Proposition
5). The $\si_j$ are those ladder epochs which are visited only
once. For $\tau_k$ to be such a special ladder epoch it is
required that after $\tau_k$ the random walk stay strictly above
its value at $\tau_k$, that is, it is required that $\cC_k$ occur.
The special ladder epochs $\si_s$ are regeneration
epochs, because they separate the path of the random walk
$\{S_n\}$ into two pieces which do not overlap (except that the
endpoint of one of these pieces coincides with the initial point
of the next piece).



\bigskip
On the event $\{\nu_j < \infty\}$ we define the $j$-th {\it
excursion} $\Up_j$ to be the sequence of random variables $(S_n -
S_{\tau_{\nu_j}}) = (S_n-\nu_j), \nu_j \le n < \nu_{j+1}$. We
already proved in Remark 2 that all $\nu_j$ are finite a.s. $[Q]$.
To describe the distribution of the excursions we introduce some
collections of possible finite sequences which can be the value of
$\Up_j$. For $1 \le m < \infty$, we define  $\wt \cM_m$ as the
collection of sequences $\boldsymbol \eta = (\eta_0,\eta_1,
\eta_2, \dots, \eta_m)$ which satisfy
\begin{equation}\label{new12}
\eta_0 = 0, \eta_i - \eta_{i-1} \in \{+1, -1\} \text{ for }1 \le i
\le m,
\end{equation}
\begin{equation}\label{new13}
\text{for any $x \in \Bbb Z, \eta_i = x$ for at most $L_0$ values
of }i \in [0,m],
\end{equation}
\begin{equation}\label{new14}
\eta_m > \eta_i \text{ for } 0 \le i < m,
\end{equation}
but there is no $0< j < m$ such that
\begin{equation}\label{new14.1}
S_i < S_j < S_\ell < S_m \text{ for } i<j < \ell< m.
\end{equation}



\bigskip
These collections will serve to describe the distribution of
$\Up_j$ when $j=0$. For $j \ge 1$ we shall use $\cM_m$ which is
defined as the collection of sequences $\boldsymbol \eta =
(\eta_0,\eta_1, \eta_2, \dots, \eta_m) \in \wt \cM_m$ which in
addition satisfy
\begin{equation}\label{new15}\eta_i > 0 \text{ for } 1 \le i \le m.
\end{equation}

\begin{prop} Under $Q$ all the $\nu_j$ are a.s. finite.
Moreover, the excursions $\Up_j$ are independent, with a distribution specified by:
\begin{equation}\label{new30}
Q(\Up_0 = \boldsymbol \eta = (\eta_0, \dots, \eta_m)) = \frac
1{\wt Z} e^{C_4m}P((S_0, \dots,S_m) = \boldsymbol \eta)
\end{equation}
for any $\boldsymbol \eta \in \wt \cM_m$.
Here $\wt Z$ is a normalizing factor given by
\begin{equation}\label{new31}
\wt Z = \sum_{m = 0}^\infty e^{C_4m} \sum_{\boldsymbol \eta \in
\wt \cM_m} P((S_0, \dots,S_m) = \boldsymbol \eta).
\end{equation}
Similarly, the distribution of $\Up_s$ with $s \ge 1$ is given by
\begin{equation}\label{new38}
Q(\Up_s = \boldsymbol \eta = (\eta_0, \dots, \eta_m)) = \frac 1Z
e^{C_4m}P((S_0, \dots, S_m) = \boldsymbol \eta)
\end{equation}
for any $\boldsymbol \eta \in \cM_m$, with $Z$ given by
\begin{equation}\label{new39}
Z = \sum_{m=0}^\infty e^{C_4m} \sum_{\boldsymbol \eta \in \cM_m}
P((S_0, \dots, S_m) = \boldsymbol \eta).
\end{equation}
In particular, the $(\Up_s, s \ge 1)$, are i.i.d. under $Q$.
Moreover, for every $s\ge 0$,
\begin{equation}\label{new25}
E^Q(\nu_{s+1} - \nu_s) < \infty.
\end{equation}
\end{prop}
\begin{proof} We already know from Remark 2 that all $\nu_j$ are
finite a.s. $[Q]$.

Now suppose that $\cH(\boldsymbol \eta)$ occurs for some $\eta \in
\wt \cM_m$. By (\ref{new14}) we then automatically have that $m =
\tau_s$ for $s = \eta_m$. Therefore, on $\cH(\boldsymbol \eta)
\cap \cC_s = \cH(\boldsymbol \eta) \cap \cC_{\eta_m}$, $\tau_s
\cap \cC_s$ occurs and $s= \eta_m$ has to equal $\si_r$ for some
$r$ and $s$ has to be one of the $\nu_j$. In fact (\ref{new14.1})
shows that there can be no $j < m$  such that $\eta_j$ is an
earlier $\si$, i.e. $\si_t$ with $t < r$. Thus, on
$\cH(\boldsymbol \eta) \cap \cC_s$ it holds that $\si_1 = \tau_s$.
Moreover, $\si_1 = \tau_s$ can occur only if $\cH(\boldsymbol
\eta)$ for some $\boldsymbol \eta \in \wt \cM_m$ occurs, as well
as $\cC_s$. Thus, $\Up_0 = \boldsymbol \eta$ is possible only if
$\boldsymbol \eta$ lies in $\wt \cM_m$ for some $m$. Furthermore
$$
\{\Up_0 = \boldsymbol \eta\}= \{\cH(\boldsymbol \eta) \cap
\cC_{\eta_m}\} \text{ for } \boldsymbol \eta \in \wt \cM_m.
$$
Also, for $\boldsymbol \eta \in \wt \cM_m$, $\cH(\boldsymbol \eta)
\in \cG_{\eta_m}$ (because $m = \tau_{\eta_m}$ by (\ref{new14})).
Hence (\ref{new6.7z}), with $\cE$ replaced by $\cH(\boldsymbol \eta)$,
shows that
$$
Q(\Up_0 = \boldsymbol \eta) = Q(\cH(\boldsymbol \eta) \cap
\cC_{\eta_m}) = C_5e^{C_4m}P((S_0, \dots, S_m) = \boldsymbol
\eta, \sup _xL(m,x) \le L_0).
$$
The condition $\sup_x L(m,x) \le L_0$ can be dropped here, because
this is automatic if $S_i = \eta_i, 0 \le i \le m$ for some
$\boldsymbol \eta \in \wt \cM_m$ (by (\ref{new13})). This implies
(\ref{new30}) with (\ref{new31}).

To prove the statements (\ref{new38}) and (\ref{new39}) in
Proposition 5 we have to show that for $\boldsymbol \eta^{(0)} \in
\wt \cM_{m(0)}$, and $\boldsymbol \eta^{(s)} \in \cM_{m(s)}$ for
$1 \le s \le r$, it holds for some constant $C$
\begin{equation}\label{new32}
Q(\Up_s = \boldsymbol \eta^{(s)}, 0 \le s \le r)= C
\prod_{s=0}^r \big[e^{C_4m(s)}P((S_0, \dots, S_{m(s)}) =
\boldsymbol \eta^{(s)})\big]
\end{equation}
Let $\boldsymbol \eta^{(s)} = (\eta^{(s)}_0 = 0, \dots,
\eta^{(s)}_{m(s)})$ and write
$$
q(s) = \sum_{j=0}^{s-1} m(j)
$$
(with $q(0) =0$). Then the event in the left hand side of
(\ref{new32}) will occur if and only if
\begin{equation}\label{new33}
\cE^{(r)}:= \Big\{S_{q(s)+i} = \sum_{j=0}^{s-1} \eta^{(j)}_{m(j)}
+ \eta^{(s)}_i,  \quad 0 \le i \le m(s), 0 \le s \le r \Big\},
\end{equation}
as well as
\begin{equation}\label{new34}
\bigcap_{s=0}^r \cC\Big(\sum_{j=0}^s \eta^{(j)}_{m(j)}\Big)
\end{equation}
occur. Here we have written $\cC(a)$ for $\cC_a$ to avoid
complicated subscripts. By the definitions of $\wt \cM_m(0)$ and
the $\cM_{m(s)}$, $\eta^{(0)}_{m(0)} > 0$ and all $\eta^{(s)}_j, 0
\le j \le m(s)$ are nonnegative. Therefore $\sum_{j=0}^s
\eta_{m(j)}^{(j)} \ge \sum_{j=0}^{s-1} \eta_{m(j)}^{(j)}, 1 \le s
\le r$, and
$$
\bigcap_{s=0}^{r-1} \cC\Big(\sum_{j=0}^s \eta^{(j)}_{m(j)}\Big)
$$
is contained in
$$
\cE^{(r)} \bigcap \cC\Big(\sum_{j=0}^r \eta^{(j)}_{m(j)}\Big).
$$
Also the event (\ref{new33}) is contained in $\cG\big(\sum_{s=0}^r
\eta^{(s)}_{m(s)}\big)$ (where we have written $\cG(a)$ for
$\cG_a$), since on (\ref{new33}) $S_.$ reaches the level
$\sum_{s=0}^r \eta^{(s)}_{m(s)}$ first at the time $q(r+1)$. It
now follows from the fact that the value of (\ref{new15.1z}) is
given by (\ref{new15zz}) that the left hand side of (\ref{new32})
equals
\begin{equation}\label{new35}
C_5e^{C_4q(r+1)}P(\cE^{(r)}, \sup_xL(q(r+1), x) \le L_0)
\end{equation}
Finally,
$$
\cE^{(r)} = \bigcap_{s=0}^r \big\{S_{q(s)+i} - S_{q(s)} =
\eta^{(s)}_i, \quad 0 \le i \le m(s)\big\},
$$
and on $\cE^{(r)}$ the range of $\{S_{q(s)+i}, 0 \le i \le m(s)\}$
consists of the integers in the interval $[\sum_{j=0}^{s-1}
\eta^{(j)}_{m(j)},  \sum_{j=0}^s \eta^{(j)}_{m(j)}], \quad 1 \le s
\le r$. The interiors of these intervals are disjoint and any
value $x$ in the interior of these intervals is taken on at most
$L_0$ times by $\{S_{q(s)+i}, 0 \le i \le m(s)\}$ if
$(S_{q(s)}-S_{q(s)}, S_{q(s)+1}-S_{q(s)},\dots, S_{q(s+1)}) -
S_{q(s)}) = \boldsymbol \eta^{(s)}$, by virtue of (\ref{new13}).
Moreover, on $\cE^{(r)}$, the endpoints $\sum_{j=0}^s
\eta^{(j)}_{m(j)}, \quad 0 \le s \le r$ are even taken on only
once by the $S_i, 0 \le i \le q(r+1)$, because $\sum_{j=0}^s
\eta^{(j-1)}_{m(j-1)} > \sum_{j=0}^{s-1} \eta^{(j)}_{m(j)}$.
Therefore, the condition $\sup_x L(q(r+1),x) \le L_0$ is
automatically fulfilled on $\cE^{(r)}$ and can be dropped from
(\ref{new35}). The result is
\begin{align}
Q(\Up_r = \boldsymbol \eta^{(r)}, 0 \le r \le s)
& =C_5e^{C_4q(r+1)}P(\cE^{(r)})\nonumber\\
& = C_5 \prod_{s=0}^r \big[e^{C_4m(s)} P(S_i = \eta^{(s)}_i, 0
\le i \le m(s))\big]
\label{new37}
\end{align}
for $\boldsymbol \eta^{(0)} \in
\wt\cM_{m(0)},\boldsymbol \eta^{(s)} \in \cM_{m(s)}, 1 \le s \le r$. The fact that the right hand side here is a product of factors
each of which depend on the value of one $\Up_s$ only shows that
the $\Up_s$ are independent. The actual distribution of the
$\Up_s$ can also be read of from (\ref{new37}) and is given by
(\ref{new38}) and (\ref{new39}).

Finally, the random variables $(\nu_{s+1} -\nu_s), \; s \ge 1$,
are i.i.d. under $Q$, so that by the renewal theorem,
\begin{equation}\label{ren}
\frac 1n \sum_{\ell =1}^n Q(\ell \text{ equals some }\nu_j)
\to [E^Q(\nu_2-\nu_1)]^{-1} \text{ as } n \to \infty.
\end{equation}
However, by (\ref{new21.2z}), we know that
$$
Q(\ell \text{ equals some }\nu_j)\ge Q(\cC_\ell) \ge C_3
$$
hence
$$
\liminf_{n\to \infty}\frac 1n \sum_{\ell =1}^n Q(\ell \text{ equals some }\nu_j) \ge C_3
$$
Since $C_3>0$, this and (\ref{ren}) imply (\ref{new25}). \qed
\end{proof}
\begin{cor}
\begin{equation}\label{new38zz}
E^Q(\si_{s+1}- \si_s) = E^Q(\tau_{\nu_{s+1}}-\tau_{\nu_s}) <
L_0 E^Q(\nu_{s+1} - \nu_s) < \infty.
\end{equation}
\end{cor}
\begin{proof} By (\ref{new11}), the amount of time spent by the walk in
any interval $[a, b) \subset \Bbb Z$ is at most $L_0(b-a)$. By
definition of the $\tau's$ and the $\si$'s the walk stays in the
interval $[\nu_j, \nu_{j+1})$ during $[\tau_{\nu_j},
\tau_{\nu_{j+1}}) = [\si_j, \si_{j+1})$. Thus (\ref{new38zz})
follows from (\ref{new25}). The strict inequality in (\ref{new38zz}) follows from the fact that for every $j\ge 0$, every site $x$ between two successive regeneration levels $x \in [\nu_j,\nu_{j+1}) \cap \mathbb{Z}$ is visited at most $L_0$ times, except $x=\nu_j$ itself which is visited at most once. It follows that $$\sigma_{j+1}-\sigma_j \le L_0(\nu_{j+1} - \nu_j -1) +1 $$ almost surely. Taking expectations leads to the strict inequality in (\ref{new38zz}). \qed
\end{proof}
With this in mind, routine manipulations show that under $Q$, the position $X_t$ satisfies the law of large numbers
\begin{equation}
\limn\frac{X_n}n := \gamma(L_0)
\end{equation}
exists almost surely under $Q$, with
$$
\gamma(L_0):= \frac{E^Q(\nu_{s+1} - \nu_s)}{E^Q(\si_{s+1}- \si_s)} .
$$
By (\ref{new38zz}), we see that
$$
\gamma(L_0)>1/L_0
$$
which concludes the proof of Theorem \ref{T:rw}.\end{proof}

\end{document}